\documentclass[11pt,a4paper]{article}

\usepackage{amsmath, amsfonts, amssymb, amsthm, mathtools, mathrsfs}
\numberwithin{equation}{section}
\usepackage[margin=1in]{geometry}

\usepackage{hyperref}
\hypersetup{hidelinks}

\usepackage{setspace}
\usepackage{needspace}

\providecommand{\keywords}[1]{%
\par\smallskip
\noindent\textit{Keywords:} #1\par
}

\providecommand{\subjclass}[1]{%
\par\smallskip
\noindent\textit{MSC (2020):} #1\par
}

\newtheorem{assumption}{Assumption}

\theoremstyle{plain}
\newtheorem{thm}{Theorem}[section]
\newtheorem{lem}[thm]{Lemma}
\newtheorem{prop}[thm]{Proposition}
\newtheorem{corollary}[thm]{Corollary}
\theoremstyle{definition}
\newtheorem{definition}[thm]{Definition}
\theoremstyle{remark}
\newtheorem{remark}{Remark}

\title{Sticky Walsh Diffusions with General Spinning Measures}
\author{Wajdi Touhami\\[0.5em]
\small University of Tunis El Manar, Higher Institute of Medical Technologies\\
\small 9 street Dr. Zouhair Essafi, 1006 Tunis, Tunisia\\
\small \texttt{wajdi.touhami@istmt.utm.tn}}
\date{}
\begin{document}
 \maketitle

\begin{abstract}
We study sticky Walsh diffusions with a general, possibly non-atomic, spinning measure and ray-dependent drift and diffusion coefficients. Stickiness is encoded by an occupation-time/local-time relation and the angular dynamics by a local-time splitting property, leading to an integral Wentzell boundary condition at the vertex.

We characterize these diffusions through a de-stickifying time change and establish existence and uniqueness in distribution under suitable assumptions. Using the generalized Freidlin--Sheu formula, we obtain the associated martingale problem and a local Dynkin formula. Under additional global regularity and uniform ellipticity assumptions, we identify the full infinitesimal-generator domain on a Banach space of uniformly tree-continuous functions vanishing at radial infinity and derive an explicit resolvent formula.

We also study small-radius exits from the vertex, obtaining first-order asymptotics for the mean exit time and expected vertex occupation time. We derive an exact joint Laplace transform of the exit time, vertex occupation time, and outgoing direction. This shows that the vertex occupation time is exponentially distributed and independent of the outgoing direction, and yields an explicit scale-function formula for the angular exit distribution. As the exit radius \(\delta\) tends to zero, this distribution converges to the spinning measure in total variation at rate \(O(\delta)\), while the normalized exit time and direction converge jointly in distribution.
\end{abstract}



\section{Introduction}
\label{sec:introduction}

Walsh Brownian motion is a natural extension of one-dimensional
reflected Brownian motion to a family of rays meeting at a common
vertex. Away from the vertex, the process evolves along a single ray,
whereas each visit to the vertex allows a new direction to be selected
according to a prescribed probability measure, usually called the
\emph{spinning measure}. This construction goes back to Walsh and was
further developed by Barlow, Pitman and Yor; see
\cite{BarlowPitmanYor1989}. It provides a basic model for diffusion
processes on ramified state spaces and is closely related to diffusion
processes on graphs.

A stochastic calculus for Walsh-type processes has subsequently been
developed in several directions. Freidlin and Sheu
\cite{FreidlinSheu2000} established an It\^o-type formula for
diffusions on finite graphs. Hajri and Touhami
\cite{HajriTouhami2014} extended the corresponding formula for Walsh
Brownian motion to an arbitrary probability measure on the angular
space, thus allowing the spinning measure to have non-finite and, in
particular, non-atomic support. Ichiba, Karatzas, Prokaj and Yan
\cite{IchibaKaratzasProkajYan2018} developed stochastic integral
equations and change-of-variable formulas for Walsh semimartingales.
A general theory of Walsh diffusions with ray-dependent drift and
diffusion coefficients was developed by Karatzas and Yan
\cite{KaratzasYan2019}, including existence and uniqueness in
distribution, local-time splitting among sets of rays, and
change-of-variable formulas. Bayraktar, Zhang and Zhang
\cite{BayraktarZhangZhang2025} constructed Walsh spider diffusions
on finite star-shaped graphs by time-changing multi-parameter
processes, starting from prescribed infinitesimal generators.
Related questions concerning the generator on a class of test
functions, stationary distributions, and convergence for Walsh
diffusions were studied by Ichiba and Sarantsev
\cite{IchibaSarantsev2019}.

More recently, Martinez and Ohavi \cite{MartinezOhavi2025}
established a well-posed martingale problem for Walsh spider
diffusions on finite star graphs with spinning probabilities and
coefficients allowed to depend on the accumulated local time at the
vertex. In a subsequent work \cite{MartinezOhavi2026}, they obtained
an It\^o formula, Markov properties, and a Feynman--Kac
representation in that setting. Their local-time-dependent scattering
mechanism is distinct from the sticky occupation mechanism studied
here: our spinning measure is fixed but may be non-atomic, while the
process can spend positive Lebesgue time at the vertex.

In the Walsh diffusion frameworks just described, the vertex is
typically taken to be instantaneous. In particular, in the framework
of \cite{KaratzasYan2019}, a non-stickiness condition is imposed in
order to exclude positive occupation time at the origin. This is in contrast
with \emph{sticky diffusions}, for which the process may spend a
positive amount of Lebesgue time at a distinguished boundary point.
Sticky behavior is classical for one-dimensional diffusions and can
be described either through an occupation-time/local-time relation or
through a suitable time change.

Sticky behavior has also been investigated on metric and star graphs.
The general theory of diffusions on graphs leads naturally to
Feller--Wentzell type conditions at the vertices; see, for instance,
Freidlin and Wentzell \cite{FreidlinWentzell1993} and Freidlin and
Sheu \cite{FreidlinSheu2000}. Bonaccorsi and D'Ovidio
\cite{BonaccorsiDOvidio2024} constructed sticky Brownian motions on
finite star graphs by slowing down the process at the vertex and
investigated related trapping mechanisms. More recently, Berry and
Colantoni \cite{BerryColantoni2026} studied sticky diffusions on
finite star graphs. They characterized such processes as time changes
of non-sticky diffusions, obtained a stochastic differential equation
for the radial component, and established a sticky version of the
Freidlin--Sheu formula. A time-change representation for more general
regular diffusions on star-shaped graphs, including sticky gluing
conditions, was also developed by Anagnostakis
\cite{Anagnostakis2025}. Bobrowski and Pilipenko
\cite{BobrowskiPilipenko2025} described Brownian motions on
star graphs with finitely many edges, including their generator
and resolvent formulations. Basti, Berry and Camilli
\cite{BastiBerryCamilli2026} developed Markov-chain approximations
for sticky diffusions on networks and numerical schemes for the
associated Hamilton--Jacobi--Bellman equations. Related
approximation results on finite star-like graphs were obtained by
Bobrowski and Ratajczyk
\cite{BobrowskiRatajczyk2026}, who studied the approximation of
Walsh spider processes by snapping-out Brownian motions.

These strands of work leave a natural question. The general Walsh
diffusion framework allows an arbitrary spinning measure and
ray-dependent coefficients, but the vertex is instantaneous, whereas
the sticky star-graph framework allows positive occupation time at the
vertex but is formulated for a finite collection of edges. It is
therefore natural to ask whether stickiness can be incorporated into
the general Walsh setting while preserving the angular splitting
mechanism associated with the spinning measure.

Passing from finitely many edges to an arbitrary spinning measure is
not simply a matter of replacing a finite weighted sum by an
integral. With infinitely many rays, the Walsh space equipped with
its tree metric is not locally compact at the vertex. Consequently,
generator-domain descriptions based on the usual $C_0$-theory for
finite metric graphs cannot be imported without qualification.
Moreover, solving a boundary problem separately on each ray does not
by itself ensure that the solutions and their boundary derivatives
are jointly Borel measurable in the radius and direction, or that
their bounds are uniform over the support of a non-atomic measure.
Such properties are needed to define and control the integral
boundary term. On the probabilistic side, the insertion or removal
of sticky time must remain compatible with the local-time splitting
identity for every Borel set of directions. These are the specific
issues that our time-change argument and uniform local estimates
address; we do not rely on an automatic extension of finite-graph
generator theory.

In this paper, we study sticky Walsh diffusions with a general spinning measure \(\nu\) on the angular space, which may in particular be non-atomic, and with drift and diffusion coefficients depending on both the radial and angular variables. Thus, the state space is not restricted to a finite
collection of edges. Instead, we work on the Walsh space
\[
\mathcal W
=
\{v\}
\cup
\bigl((0,\infty)\times S^1\bigr),
\]
equipped with its natural tree topology. A point different from the
vertex is represented by $(r,\theta)$, where $r>0$ denotes the radial
coordinate and $\theta\in S^1$ the direction.

Let
\(Y=(R,\Theta)\)
denote the process. Away from the vertex, its radial component evolves
according to the ray-dependent operator
\[
\mathcal{L}_\theta f_\theta(r)
=
b_\theta(r) f_\theta'(r)
+
\frac12 \sigma_\theta^2(r) f_\theta''(r).
\]
The angular variable remains constant during each excursion away from
the vertex. The distribution of excursions among sets of directions
is prescribed by a probability measure $\nu$, through the angular
local-time splitting property
\[
L_t^0(R_\Lambda)
=
\nu(\Lambda)L_t^0(R),
\qquad
\Lambda\in\mathcal B(S^1),
\]
in the appropriate localized sense. Stickiness at the vertex is
described by the occupation identity
\[
\int_0^t
\mathbf 1_{\{Y_s=v\}}\,ds
=
\eta L_t^0(R),
\qquad
\eta>0,
\]
where $L^0(R)$ denotes the local time of the radial component at
the origin and $\eta$ is the stickiness parameter. The limiting
case $\eta=0$ corresponds to the non-sticky regime, whereas the
present paper considers $\eta>0$.

The combination of angular local-time splitting according to $\nu$
and sticky occupation at the vertex leads naturally to the integral
Wentzell boundary condition
\[
\eta\,\mathcal{L} f(v)
=
\int_0^{2\pi}
f_\theta'(0+)\,\nu(d\theta).
\]
For a spinning measure with finite support, the integral reduces to
the familiar weighted sum of outgoing derivatives occurring in
diffusions on finite star graphs. For a general measure $\nu$, the
boundary condition simultaneously describes the sticky delay at the
vertex and the distribution of local time among the possible
directions.

One of the main results of the paper is a time-change
characterization of the process. By deleting the Lebesgue time spent
at the vertex, we obtain a non-sticky Walsh diffusion with the same
ray-dependent coefficients and the same spinning measure. The
resulting process falls within the general Walsh diffusion framework
of Karatzas and Yan \cite{KaratzasYan2019}. Conversely, the sticky
process is recovered from the corresponding non-sticky diffusion by
inserting an amount of physical time proportional to its local time
at the vertex. This provides, in particular, uniqueness in
distribution for the sticky process under the corresponding
assumptions on the coefficients.

The de-sticky representation also allows us to use the
Freidlin--Sheu change-of-variable formula for the underlying
non-sticky Walsh diffusion to derive the martingale problem
associated with the Wentzell operator. Under the sticky time change,
the local-time contribution in the Freidlin--Sheu formula is
converted, through the Wentzell boundary condition and the sticky
occupation identity, into the contribution of the time spent at the
vertex. This yields the corresponding Dynkin formula, which provides
the main analytic tool for the local analysis near the vertex.

Under additional, globally bounded, uniformly elliptic and radially
Lipschitz coefficients, we strengthen this analytic description:
Theorem~\ref{gen:gen} identifies the complete domain of the
infinitesimal generator on a Banach space of uniformly continuous
functions for the tree metric, subject to uniform radial vanishing.
The result gives an explicit resolvent in terms of raywise killed
resolvents and an angularly averaged Wentzell boundary term. No
continuity of the coefficients in the angular variable is imposed.
For non-sticky Walsh Brownian motion with a general spinning measure,
Chen and Fukushima \cite[Theorem~5.4]{ChenFukushima2015} characterized
the $L^2$-generator domain through an averaged zero-flux condition at
the vertex. Our result instead concerns the generator on $\mathscr B$
for sticky Walsh diffusions with ray-dependent coefficients, with an
integral Wentzell boundary condition.
This complementary result is separate from the martingale formulation
under the more general hypotheses of our construction.

Finally, we investigate the local
behavior of the process near the vertex. If
\[
\tau_\delta
:=
\inf\{t\geq0:R_t=\delta\},
\]
for a process started from $v$, we prove
\[
\mathbb E_v[\tau_\delta]
=
\eta\delta+O(\delta^2),
\qquad
\delta\downarrow0.
\]
We also establish the corresponding occupation-time asymptotic
\[
\mathbb E_v
\left[
\int_0^{\tau_\delta}
\mathbf 1_{\{Y_s=v\}}\,ds
\right]
=
\eta\delta+O(\delta^2).
\]
These estimates quantify the sticky delay near the vertex.

A further quantitative result is the exact joint Laplace transform
in Proposition~\ref{prop:joint-exit-transform}. Writing
$\Gamma_\delta:=\int_0^{\tau_\delta}\mathbf1_{\{Y_s=v\}}\,ds$,
we show that, for $\lambda,q\ge0$ and every bounded Borel function
$g$ on the angular space,
\[
\mathbb E_v\!\left[
 e^{-\lambda\tau_\delta-q\Gamma_\delta}
 g(\Theta_{\tau_\delta})\right]
=
\frac{\displaystyle\int_0^{2\pi}
 g(\theta)\psi_{\lambda,\theta}(\delta)^{-1}\,\nu(d\theta)}
{\displaystyle\int_0^{2\pi}
 \phi_{\lambda,\theta}(\delta)
 \psi_{\lambda,\theta}(\delta)^{-1}\,\nu(d\theta)
 +\eta(\lambda+q)},
\]
where $\phi_{\lambda,\theta}$ and $\psi_{\lambda,\theta}$ are
normalized fundamental solutions of the raywise equation
$\mathcal{L}_\theta w=\lambda w$. Thus the dependence on the
raywise coefficients is retained in an explicit angular average,
rather than disappearing in the small-radius asymptotics.
Setting $\lambda=0$ yields an exact exponential law for
$\Gamma_\delta$ and its independence from the outgoing direction.
Writing $s_\theta$ for the normalized raywise scale function and
$K_\delta:=\int_0^{2\pi}s_\theta(\delta)^{-1}\,\nu(d\theta)$,
the outgoing direction has the exact law
\[
\mathbb P_v(\Theta_{\tau_\delta}\in d\theta)
=\frac{s_\theta(\delta)^{-1}}{K_\delta}\,\nu(d\theta).
\]
In particular, this law converges to $\nu$ in total variation
at rate $O(\delta)$ as $\delta\downarrow0$.
The transform also leads to a joint small-radius limit for
$(\tau_\delta/\delta,\Theta_{\tau_\delta})$: the limit consists of
a scaled exponential variable and an independent direction of law
$\nu$. No independence between the \emph{total} exit time and the
outgoing direction is asserted at a fixed radius when the
coefficients depend on the ray.

To the best of our knowledge, the existing literature does not treat
simultaneously a sticky vertex, a general---possibly non-atomic---
spinning measure, and ray-dependent diffusion coefficients in the
Walsh setting considered here. Our contribution is neither the
time-change principle nor the Walsh angular mechanism in isolation,
but their rigorous compatibility in the continuum-ray setting, the
resulting integral Wentzell martingale problem, its complete
generator-domain refinement under stronger regularity assumptions,
and the quantitative exit and occupation-time identities above. In particular, the exact
joint transform complements the first-order exit estimates by
recording how the ray-dependent dynamics, spinning measure, and
sticky parameter enter the same distributional identity.

The paper is organized as follows.
Section~\ref{sec:probabilistic-construction} introduces the Walsh
space and the definition of sticky Walsh diffusion, and establishes
the de-sticky time-change characterization and uniqueness in
distribution.
Section~\ref{sec:wentzell-martingale} introduces the Wentzell
operator, derives the martingale problem and a local Dynkin formula,
and identifies the full generator domain under additional global regularity.
Section~\ref{sec:local-properties} first establishes the small-radius
asymptotics of the mean exit time from the vertex. It then derives
an exact joint Laplace transform of the exit time, vertex occupation
time, and outgoing direction, yielding an exponential occupation-time
law, an explicit angular exit distribution with total-variation
convergence to the spinning measure, and a joint small-radius limit.

\section{Probabilistic construction}
\label{sec:probabilistic-construction}

Throughout, directions in $S^1$ are represented by angles $\theta\in[0,2\pi)$, taken modulo $2\pi$. Accordingly, integrals with respect to the spinning measure $\nu$ are written with limits $0$ and $2\pi$.

Throughout the paper, we use the convention $\inf\varnothing=\infty$.

The Walsh space $\mathcal W$ introduced in the
Introduction can equivalently be represented as
\[
\mathcal W
:=
\bigl([0,\infty)\times[0,2\pi)\bigr)/{\sim},
\]
where
\((0,\theta)\sim(0,\theta')\)
for all
\(\theta,\theta'\in[0,2\pi).\) Thus all points with radial coordinate zero are identified with a
single point, denoted by \(v\), called the \emph{vertex}. This is the
usual state-space representation for Walsh-type processes; see, for
instance, \cite{BarlowPitmanYor1989,KaratzasYan2019}.

Every point \(x\in\mathcal W\setminus\{v\}\) can be represented
uniquely as
\(x=(r,\theta),
\;
r>0,\;\theta\in[0,2\pi).\)
We write
\(\|x\|=r,
\,
\arg(x)=\theta,\)
and set \(\|v\|=0\). We equip \(\mathcal W\) with the tree metric
\[
\rho\bigl((r,\theta),(s,\varphi)\bigr)
=
\begin{cases}
|r-s|,
& \theta=\varphi,\\[1mm]
r+s,
& \theta\neq\varphi,
\end{cases}
\]
with
\(\rho((r,\theta),v)=r.\)
The topology induced by \(\rho\) will be referred to as the
\emph{tree topology}; see \cite{KaratzasYan2019}.

In particular, continuity on $\mathcal W\setminus\{v\}$ is equivalent
to continuity along each ray, whereas continuity at the vertex means
\[
\sup_{\theta\in[0,2\pi)}
|f(r,\theta)-f(v)|
\longrightarrow0,
\qquad r\downarrow0.
\]
All continuity statements on \(\mathcal W\) are understood with
respect to the tree topology.

For measurability, however, we equip $\mathcal W$ with the usual
Euclidean Borel $\sigma$-field through its identification with
$\mathbb R^2$ in polar coordinates. The tree topology is used for
path continuity, while adaptedness, measurable functions, and
probability laws refer to this Borel structure.

Throughout the paper, for every nonnegative
continuous semimartingale \(Z\), the local
time \(L^0(Z)\) at zero is normalized by
Tanaka's formula
\[
Z_t
=
Z_0+
\int_0^t
\mathbf 1_{\{Z_s>0\}}\,dZ_s
+
L_t^0(Z).
\]
This convention agrees with that of
Karatzas and Yan~\cite{KaratzasYan2019}.

Motivated by the semimartingale framework for Walsh diffusions of
Karatzas and Yan \cite{KaratzasYan2019} and by the occupation-time
description of sticky diffusions
\cite{EngelbertPeskir2014,BerryColantoni2026}, we use the following
probabilistic definition.

\begin{definition}[Sticky Walsh diffusion]
\label{def:sticky-Walsh-diffusion}

Let \(b\) and \(\sigma\) be coefficient
functions such that the mappings
\[
(r,\theta)\longmapsto b_\theta(r),
\qquad
(r,\theta)\longmapsto \sigma_\theta(r)
\]
are jointly Borel measurable on
\((0,\infty)\times[0,2\pi)\). Let \(\nu\) be a probability measure on
\([0,2\pi)\), and let \(\eta>0\).

A sticky Walsh diffusion associated with \((b,\sigma,\nu,\eta)\),
defined up to an explosion time, consists of a filtered probability
space
\[
(\Omega,\mathcal F,\mathbb F,\mathbb P),
\qquad
\mathbb F=(\mathcal F_t)_{t\ge0},
\]
a one-dimensional process \(W^Y\), and a process \(Y\)
defined on the stochastic interval \([0,\tau^Y)\), where
\(\tau^Y\) is an \(\mathbb F\)-stopping time representing
its lifetime.

The process \(Y\) takes values in \(\mathcal W\) before
its lifetime. Its radial component is defined by
\(R_t^Y:=\|Y_t\|,\; 0\le t<\tau^Y.\)
For each \(n\ge1\), define
\[
\tau_n^Y:=\inf\{0\le t<\tau^Y:R_t^Y\ge n\}.
\]
We assume that
\(\tau^Y=\lim_{n\to\infty}\tau_n^Y,\)
so that the lifetime coincides with the radial explosion time.
The process is considered only before explosion, with an optional
cemetery-state extension thereafter.

The following conditions hold:

\begin{enumerate}

\item
The filtration \(\mathbb F\) satisfies the usual conditions, and
\(W^Y\) is a standard one-dimensional Brownian motion with respect to
\(\mathbb F\).

\item
The process $Y$ is adapted to $\mathbb F$ and has continuous
paths in the tree topology on $[0,\tau^Y)$.

The angular component is constant during each excursion away
from the vertex. More precisely, on every connected component
$I$ of
\(\{t\in[0,\tau^Y):R_t^Y>0\},\)
there exists $\theta\in[0,2\pi)$ such that
\(Y_t=(R_t^Y,\theta),\; t\in I.\)

\item
For every \(n\in\mathbb N\) and \(T\ge0\),
\[
\int_0^{T\wedge \tau_n^Y}
\mathbf 1_{\{Y_t\neq v\}}
\bigl(|b(Y_t)|+\sigma^2(Y_t)\bigr)\,dt
<\infty,
\qquad\text{a.s.}
\]

\item
For every \(n\in\mathbb N\), the stopped radial process
\(R^Y_{\cdot\wedge \tau_n^Y}\) is a semimartingale satisfying
\[
\begin{aligned}
R_{T\wedge \tau_n^Y}^Y
={}&R_0^Y
+
\int_0^{T\wedge \tau_n^Y}
\mathbf 1_{\{Y_t\neq v\}}b(Y_t)\,dt\\
&+
\int_0^{T\wedge \tau_n^Y}
\mathbf 1_{\{Y_t\neq v\}}\sigma(Y_t)\,dW_t^Y
+
L_T^0\!\left(R^Y(\cdot\wedge \tau_n^Y)\right),
\end{aligned}
\]
where \(L^0(R^Y(\cdot\wedge \tau_n^Y))\) denotes the local time at
\(0\) of the stopped radial process.

\item
For every \(\Lambda\in\mathcal B([0,2\pi))\), define
\[
R_\Lambda^Y(t)
:=
\begin{cases}
R_t^Y\mathbf 1_\Lambda(\arg(Y_t)),
&Y_t\neq v,\\[1mm]
0,&Y_t=v.
\end{cases}
\]
For every \(n\in\mathbb N\), the stopped process
\(R_\Lambda^Y(\cdot\wedge\tau_n^Y)\)
is a continuous semimartingale.

\item
For every \(n\in\mathbb N\) and \(T\ge0\),
the sticky occupation identity
\begin{equation}\label{eq:sticky-occupation}
\int_0^{T\wedge\tau_n^Y}
\mathbf 1_{\{Y_t=v\}}\,dt
=
\eta
L_T^0\!\left(R^Y(\cdot\wedge\tau_n^Y)\right)
\end{equation}
holds almost surely.

\item
For every \(\Lambda\in\mathcal B([0,2\pi))\),
\(n\in\mathbb N\), and \(T\ge0\),
the local-time partition identity
\[
L_T^0\!\left(R_\Lambda^Y(\cdot\wedge\tau_n^Y)\right)
=
\nu(\Lambda)
L_T^0\!\left(R^Y(\cdot\wedge\tau_n^Y)\right)
\]
holds almost surely.

\end{enumerate}
\end{definition}

\medskip

All identities involving stopped processes are
understood up to the corresponding radial
localization times. If the explosion time is
finite, the process is considered only before
explosion, with an optional cemetery-state
extension thereafter.

The definition combines the local-time splitting property of Walsh
semimartingales, as in
\cite{IchibaKaratzasProkajYan2018,KaratzasYan2019}, with the sticky
occupation mechanism at the vertex; compare
\cite{BerryColantoni2026}.

Throughout this paper, the term non-sticky Walsh diffusion
refers to a Walsh diffusion in the framework of
Karatzas and Yan \cite{KaratzasYan2019}, with radial
coefficients $b$ and $\sigma$ and spinning measure $\nu$.

We recall that such a diffusion spends zero Lebesgue
time at the vertex. We shall use the corresponding
semimartingale and local-time properties in the
time-change characterization below.

The following change-of-variables identity is based on the standard
time-change mechanism for sticky diffusions; see, for instance,
\cite{EngelbertPeskir2014,BerryColantoni2026}.

\begin{lem}
\label{lem:change-of-variables}

Let $X$ be a non-sticky Walsh diffusion
on $\mathcal W$, with lifetime $\tau^X$,
and write
\(R_t:=\|X_t\|.\)
For $\eta>0$, define the sticky clock
\begin{equation}\label{eq:sticky-clock}
A_t:=t+\eta L_t^0(R),
\qquad 0\le t<\tau^X,
\end{equation}
and set
\(\zeta_A:=\lim_{t\uparrow\tau^X}A_t.\)
For $0\le t<\zeta_A$, let
\[
T_t:=\inf\{s<\tau^X:A_s>t\},
\qquad
Y_t:=X_{T_t}.
\]
Then, for every nonnegative Borel function
\(h:\mathcal W\to[0,\infty]\),
\begin{equation}\label{eq:sticky-change-of-variables}
\int_0^{T_t} h(X_s)\,ds
=
\int_0^t
h(Y_u)\mathbf 1_{\{Y_u\neq v\}}\,du,
\qquad 0\le t<\zeta_A.
\end{equation}
Moreover,
\begin{equation}\label{eq:sticky-vertex-occupation}
\int_0^t
\mathbf 1_{\{Y_u=v\}}\,du
=
\eta L_{T_t}^0(R),
\qquad 0\le t<\zeta_A.
\end{equation}

\end{lem}

\begin{proof}

Since $A$ is continuous and strictly
increasing on $[0,\tau^X)$, its inverse
$T$ is continuous and strictly increasing
on $[0,\zeta_A)$.

By the change-of-variables formula for
Lebesgue--Stieltjes integrals, for every
nonnegative Borel function $h$ and every
$0\le t<\zeta_A$,
\[
\int_0^t h(Y_u)\,du
=
\int_0^{T_t}h(X_s)\,dA_s.
\]

Therefore,
\[
\int_0^t h(Y_u)\,du
=
\int_0^{T_t}h(X_s)\,ds
+
\eta
\int_0^{T_t}h(X_s)\,dL_s^0(R).
\]

Since $dL_s^0(R)$ is supported on
$\{s:R_s=0\}=\{s:X_s=v\}$,
taking $h=\mathbf 1_{\{v\}}$ in the
change-of-variables formula and using
the non-stickiness of $X$, we obtain
\[
\int_0^t
\mathbf 1_{\{Y_u=v\}}\,du
=
\eta L_{T_t}^0(R).
\]

Now let $h:\mathcal W\to[0,\infty]$
be any nonnegative Borel function.
Define
\[
g(x):=
\begin{cases}
h(x), & x\neq v,\\
0, & x=v.
\end{cases}
\]
Applying the change-of-variables formula
to $g$, we obtain
\[
\int_0^t
h(Y_u)\mathbf 1_{\{Y_u\neq v\}}\,du
=
\int_0^{T_t}g(X_s)\,dA_s\\
=
\int_0^{T_t}g(X_s)\,ds\\
=
\int_0^{T_t}h(X_s)\,ds.
\]
Indeed, the second equality follows
from the fact that $g(v)=0$ and
$dL_s^0(R)$ is supported on
$\{s:X_s=v\}$, while the last equality
follows from the non-stickiness of $X$.

\end{proof}

\begin{lem}
\label{lem:desticky-clock}

Let \(Y\) be a sticky Walsh diffusion with explosion
time \(\tau^Y\), and define
\[
C_t
=
\int_0^t \mathbf 1_{\{Y_s\neq v\}}\,ds,
\qquad 0\le t<\tau^Y.
\]
Set
\(\zeta_C
=
\lim_{t\uparrow\tau^Y}C_t.\) Then \(C\) is almost surely continuous and
strictly increasing on \([0,\tau^Y)\).

\noindent{Consequently} its inverse
\[
S_t
=
\inf\{u<\tau^Y:C_u>t\},
\qquad 0\le t<\zeta_C,
\]
is continuous and strictly increasing, satisfies
\[
S_0=0,
\qquad
C_{S_t}=t,
\quad 0\le t<\zeta_C,
\]
and
\[
S_{C_t}=t,
\quad 0\le t<\tau^Y.
\]

\noindent{Moreover}, for each deterministic $t\ge0$,
the inverse $S_t$, extended by the value
$\infty$ on $\{t\ge\zeta_C\}$, is an
extended-valued stopping time. On the event $\{t<\zeta_C\}$,
the inverse $S_t$ is finite.

\end{lem}

\begin{proof}
Continuity and monotonicity of \(C\) are immediate. Suppose that \(C_b=C_a\) for some
\(0\le a<b<\tau^Y\). Then
\[
\int_a^b\mathbf 1_{\{Y_s\neq v\}}\,ds=0.
\]
By continuity of \(Y\), this implies
\(Y_s=v\) throughout \([a,b]\).

Hence \(R^Y=\|Y\|\) is constant on
\([a,b]\), so its local time at zero
does not increase there. The sticky
occupation identity \eqref{eq:sticky-occupation} therefore yields
\[
b-a
=
\eta\bigl(
L_b^0(R^Y)-L_a^0(R^Y)
\bigr)
=
0,
\]
a contradiction. Thus \(C\) is strictly increasing. Since \(C\) is continuous and strictly
increasing on \([0,\tau^Y)\), its inverse
\(S\) is continuous and strictly increasing
on \([0,\zeta_C)\), with
\[
C_{S_t}=t,
\qquad 0\le t<\zeta_C,
\]
and
\[
S_{C_t}=t,
\qquad 0\le t<\tau^Y.
\]
Moreover, for every $t<\zeta_C$,
the inverse time $S_t$ is finite.

To verify the stopping-time property,
extend $C$ to $[0,\infty)$ by setting
$C_u:=\zeta_C$ for $u\ge\tau^Y$
on $\{\tau^Y<\infty\}$.
The extended process is continuous,
adapted, and nondecreasing.

For each deterministic $t\ge0$, define
\[
S_t:=\inf\{u\ge0:C_u>t\}.
\]
For every $u\ge0$,
\[
\{S_t\le u\}
=
\bigcap_{n\ge1}\{C_{u+1/n}>t\}
\in\mathcal F_u,
\]
by the right-continuity of the filtration.
Thus $S_t$ is an extended-valued
stopping time. On $\{t<\zeta_C\}$,
it is finite and coincides with
the inverse defined in the statement.

All subsequent time changes are
understood locally before $\zeta_C$.

\end{proof}

\begin{lem}
\label{lem:excursion-preservation}

Let $Y$ be a sticky Walsh diffusion with vertex
$v$ and explosion time $\tau^Y$.
Let $C$, $S$, and $\zeta_C$ be as in
Lemma~\ref{lem:desticky-clock}, and set
\[
X_t:=Y_{S_t},
\qquad 0\le t<\zeta_C.
\]
Then the excursions of $X$ and $Y$ away
from the vertex are in one-to-one
correspondence. More precisely, if $(g,d)$
is an excursion interval of $Y$, then
$(C_g,C_d)$ is the corresponding excursion
interval of $X$, and
\[
X_{C_g+s}=Y_{g+s},
\qquad 0<s<d-g.
\]
Here, if $d=\tau^Y$, we use the convention
$C_d:=\zeta_C$. In particular, corresponding excursions
have the same duration and angular mark.
\end{lem}

\begin{proof}
Let \((g,d)\) be an excursion interval of \(Y\). Since \(Y_s\neq v\)
for \(g<s<d\),
\[
C_u-C_g
=
\int_g^u \mathbf 1_{\{Y_s\neq v\}}\,ds
=
u-g,
\qquad g\le u\le d.
\]
Hence
\(C_u=C_g+u-g\) and \(C_d-C_g=d-g.\) For \(t\in(C_g,C_d)\), set
\(u:=g+t-C_g.\)
Then \(u\in(g,d)\) and \(C_u=t\). By the definition of \(S\),
\(S_t=u=g+t-C_g,\)
and therefore
\[
X_t=Y_{g+t-C_g}\neq v.
\]
Thus \((C_g,C_d)\) is an excursion interval of \(X\). Moreover, taking
\(t=C_g+s\) gives
\[
X_{C_g+s}=Y_{g+s},
\qquad 0<s<d-g.
\]
Since \(C_d-C_g=d-g\), the corresponding excursions have the same
lifetime and angular mark.

Conversely, let \((a,b)\) be an excursion interval of \(X\), and choose
\(t_0\in(a,b)\). Since
\(X_{t_0}=Y_{S_{t_0}}\neq v,\)
the time \(S_{t_0}\) belongs to a unique excursion interval \((g,d)\)
of \(Y\). By the preceding argument, \((C_g,C_d)\) is an excursion
interval of \(X\) containing \(t_0\). Hence
\((a,b)=(C_g,C_d).\)
This proves the one-to-one correspondence.
\end{proof}

The partition of local time according to the spinning measure is a
basic feature of Walsh semimartingales; see
\cite{IchibaKaratzasProkajYan2018,KaratzasYan2019}. We next show that
this property is preserved by the de-stickifying time change.

\begin{lem}
\label{lem:local-time-partition-preservation}

Let $Y$ be a sticky Walsh diffusion with spinning
measure $\nu$, and set
\[
X_t:=Y_{S_t},
\qquad 0\le t<\zeta_C.
\]
Write $R^X=\|X\|$.
For $\Lambda\in\mathcal B([0,2\pi))$, let $R_\Lambda^Y$
be as in Definition~\ref{def:sticky-Walsh-diffusion},
and define $R_\Lambda^X$ analogously. Then, for every $0\le t<\zeta_C$,
\[
L_t^0(R^X)=L_{S_t}^0(R^Y),
\qquad
L_t^0(R_\Lambda^X)=L_{S_t}^0(R_\Lambda^Y).
\]
Consequently,
\[
L_t^0(R_\Lambda^X)
=
\nu(\Lambda)L_t^0(R^X),
\qquad \Lambda\in\mathcal B([0,2\pi)),
\]
almost surely.
\end{lem}

\begin{proof}

By definition, for every
$\Lambda\in\mathcal B([0,2\pi))$,
\[
R_t^X=R_{S_t}^Y,
\qquad
R_\Lambda^X(t)=R_\Lambda^Y(S_t),
\qquad 0\le t<\zeta_C.
\]

Let $Z$ denote either $R^Y$ or $R_\Lambda^Y$.
Since $Z$ is a continuous real-valued
semimartingale locally before explosion
and $S$ is a continuous time change
consisting of stopping times, the standard
time-change property of semimartingale
local times yields, after localization,
\[
L_t^0(Z_S)=L_{S_t}^0(Z),
\qquad 0\le t<\zeta_C.
\]
Applying this identity to $Z=R^Y$ and
$Z=R_\Lambda^Y$, we obtain
\[
L_t^0(R^X)=L_{S_t}^0(R^Y),
\qquad
L_t^0(R_\Lambda^X)=L_{S_t}^0(R_\Lambda^Y).
\]

Fix $\Lambda\in\mathcal B([0,2\pi))$.
For every $n\in\mathbb N$, the local-time
partition property in
Definition~\ref{def:sticky-Walsh-diffusion}
gives, for every $T\ge0$,
\[
L_T^0\!\left(R_\Lambda^Y(\cdot\wedge\tau_n^Y)\right)
=
\nu(\Lambda)
L_T^0\!\left(R^Y(\cdot\wedge\tau_n^Y)\right)
\]
almost surely.
By applying this identity first at rational
times and then using the continuity of
both local-time processes, we obtain an
event $\Omega_\Lambda$ of probability one on
which the identity holds simultaneously
for all $T\ge0$ and all $n\in\mathbb N$.

Since $S_t<\tau^Y$ for $t<\zeta_C$
and $\tau_n^Y\uparrow\tau^Y$, we may
evaluate the partition identity at
$T=S_t$, after localization, to obtain
\[
L_{S_t}^0(R_\Lambda^Y)
=
\nu(\Lambda)L_{S_t}^0(R^Y),
\qquad 0\le t<\zeta_C.
\]
Combining this with the time-change
identities yields
\[
L_t^0(R_\Lambda^X)
=
\nu(\Lambda)L_t^0(R^X),
\qquad 0\le t<\zeta_C.
\]
Since $\Lambda$ was arbitrary, the local-time
partition identity holds almost surely
for each fixed
$\Lambda\in\mathcal B([0,2\pi))$,
simultaneously for all
$0\le t<\zeta_C$.

\end{proof}

\begin{lem}
\label{lem:inverse-clock-localization}

Let $Y$ be a sticky Walsh diffusion, and set
\[
X_t=Y_{S_t},
\qquad
\mathcal G_t=\mathcal F_{S_t},
\qquad 0\le t<\zeta_C.
\]
For $0\le t<\tau^Y$, define
\[
M_t^Y:=
\int_0^t
\mathbf 1_{\{Y_s\neq v\}}
\sigma(Y_s)\,dW_s^Y.
\]
Then
\[
M_t^X:=M_{S_t}^Y,
\qquad 0\le t<\zeta_C,
\]
is a continuous local martingale with respect to
$(\mathcal G_t)$ on $[0,\zeta_C)$. Moreover,
\[
\langle M^X\rangle_t
=
\int_0^t\sigma^2(X_s)\,ds,
\qquad 0\le t<\zeta_C,
\]
where $\sigma(v)=0$ by convention.
\end{lem}

\begin{proof}
Define, locally before $\tau^Y$,
\(N_t:=\int_0^t\mathbf 1_{\{Y_s\ne v\}}\,dW_s^Y.\)
Then $N$ is a continuous local martingale with
\[
\langle N\rangle_t
=\int_0^t\mathbf 1_{\{Y_s\ne v\}}\,ds=C_t.
\]
By Lemma~\ref{lem:desticky-clock}, $C$ is continuous and strictly
increasing, and $S$ is its continuous inverse before $\zeta_C$.
The Dambis--Dubins--Schwarz time-change theorem
\cite[Chapter~V, \S1]{RevuzYor1999}, applied locally, therefore shows
that
\(B_t:=N_{S_t}\), $0\le t<\zeta_C$,
is a Brownian motion up to its lifetime with respect to
$\mathcal G_t=\mathcal F_{S_t}$.

For clarity about a possible explosion, let
\(U_n:=\inf\{s<\tau^Y:R_s^Y\ge R_0^Y+n\}\wedge n\),
and put
\(\gamma_n:=C_{U_n}\).
The radial-explosion convention gives
\(U_n<\tau^Y\), $U_n\uparrow\tau^Y$, and hence
$\gamma_n<\zeta_C$, $\gamma_n\uparrow\zeta_C$.
Moreover,
\(\langle N\rangle_{U_n}=C_{U_n}\le U_n\le n\).
Thus the standard stopped version of the time-change theorem applies
before each $\gamma_n$; if $\zeta_C$ is finite, the Brownian motion
may be continued on an extension of the probability space after that
time. No assertion beyond the lifetime is needed here.

With the convention $\sigma(v)=0$, associativity of stochastic
integration gives $dM_t^Y=\sigma(Y_t)\,dN_t$.
The substitution formula under the continuous time change $S$,
again applied locally, yields
\[
M_t^X=M_{S_t}^Y
=\int_0^t\sigma(Y_{S_u})\,dN_{S_u}
=\int_0^t\sigma(X_u)\,dB_u,
\qquad 0\le t<\zeta_C.
\]
Therefore $M^X$ is a continuous local $(\mathcal G_t)$-martingale,
and
\[
\langle M^X\rangle_t
=\int_0^t\sigma^2(X_u)\,du,
\qquad 0\le t<\zeta_C.
\]
\end{proof}

Throughout the characterization theorem,
non-sticky Walsh diffusions are understood
in the semimartingale framework of
Karatzas and Yan~\cite{KaratzasYan2019},
with the same normalization of local time
and the same radial coefficients and
spinning measure.

The construction from a non-sticky Walsh
diffusion presupposes the existence of such
a process. The additional condition
\textup{(LB)} is imposed only in the de-stickifying
direction to identify the explosion time.

\begin{thm}[Characterization via the de-stickifying time change]
\label{thm:characterization}
The following assertions hold
\begin{enumerate}
\item[(i)]
Assume that the coefficients satisfy the following
local lower-bound condition:
\[
\tag{LB}\label{ass:LB}
\forall n\geq1,\quad
\exists K_n\geq0
\quad\text{such that}\quad
b_\theta(r)+\frac12\sigma_\theta^2(r)\geq-K_n,
\qquad
0<r\leq n,\quad
\theta\in[0,2\pi).
\]
Let $Y$ be a sticky Walsh diffusion associated
with $(b,\sigma,\nu,\eta)$.
Let $C$, $\zeta_C$, and $S$ be as in
Lemma~\ref{lem:desticky-clock}. Then
\[
X_t:=Y_{S_t},
\qquad 0\leq t<\zeta_C,
\]
is a non-sticky Walsh diffusion associated
with $(b,\sigma,\nu)$, with explosion time
$\tau^X=\zeta_C$. Moreover,
\[
S_t=t+\eta L_t^0(\|X\|),
\qquad 0\leq t<\zeta_C.
\]

\item[(ii)]
Let $X$ be a non-sticky Walsh diffusion
associated with $(b,\sigma,\nu)$,
with explosion time $\tau^X$. Define
\[
A_s:=s+\eta L_s^0(\|X\|),
\qquad 0\leq s<\tau^X,
\]
and let
\[
\zeta_A:=\lim_{s\uparrow\tau^X}A_s,
\qquad
T_t:=\inf\{s<\tau^X:A_s>t\},
\quad 0\leq t<\zeta_A.
\]
Then
\[
Y_t:=X_{T_t},
\qquad 0\leq t<\zeta_A,
\]
is a sticky Walsh diffusion associated
with $(b,\sigma,\nu,\eta)$, with explosion
time $\tau^Y=\zeta_A$.

\end{enumerate}
\end{thm}

\begin{proof}
\noindent\textit{Assertion (i).}

\smallskip
\noindent\textit{Step 1. De-stickification.}

Let $Y$ be a sticky Walsh diffusion associated with
$(b,\sigma,\nu,\eta)$, with explosion time $\tau^Y$.
By Lemma~\ref{lem:desticky-clock}, $C$ is
continuous and strictly increasing before
explosion, and its inverse $S$ is continuous
and strictly increasing on $[0,\zeta_C)$. Moreover, each $S_t$ is a finite stopping
time and
\[
C_{S_t}=t,
\qquad 0\le t<\zeta_C.
\]
We equip $X$ with the time-changed filtration
\(\mathcal G_t:=\mathcal F_{S_t},
\; 0\le t<\zeta_C,\)
understood locally before $\zeta_C$.
Then $X$ is adapted and continuous. 
Since $X_t=Y_{S_t}$, we have
\(
R_t^X=R_{S_t}^Y,
\, 0\le t<\zeta_C
\). The radial semimartingale equation for $Y$,
interpreted locally before $\tau^Y$, is
\[
R_t^Y
=
R_0^Y
+
\int_0^t
\mathbf 1_{\{Y_u\neq v\}}b(Y_u)\,du
+
M_t^Y
+
L_t^0(R^Y),
\]
where
\[
M_t^Y
=
\int_0^t
\mathbf 1_{\{Y_u\neq v\}}
\sigma(Y_u)\,dW_u^Y.
\]
By Lemma~\ref{lem:inverse-clock-localization},
the process
\(
M_t^X:=M_{S_t}^Y,
\, 0\le t<\zeta_C,
\)
is a continuous local martingale with
respect to $(\mathcal G_t)$, with
quadratic variation
\[
\langle M^X\rangle_t
=
\int_0^t\sigma^2(X_u)\,du,
\qquad 0\le t<\zeta_C.
\]
Moreover, the proof of that lemma provides
a Brownian motion $B$ with respect to
$(\mathcal G_t)$, up to $\zeta_C$, and the
corresponding stochastic integral representation.
Taking $W^X=B$ before the lifetime
(and extending it thereafter if necessary),
we obtain
\[
M_t^X
=
\int_0^t
\mathbf 1_{\{X_u\neq v\}}
\sigma(X_u)\,dW_u^X.
\]
This representation does not require
$\sigma$ to be strictly positive. By the change-of-variables formula,
\[
\int_0^{S_t}
\mathbf 1_{\{Y_u\neq v\}}b(Y_u)\,du
=
\int_0^t
\mathbf 1_{\{X_u\neq v\}}b(X_u)\,du.
\]
Moreover, by
Lemma~\ref{lem:local-time-partition-preservation},
\(
L_{S_t}^0(R^Y)=L_t^0(R^X).
\)
Combining these identities yields
\[
R_t^X
=R_0^X
+\int_0^t
\mathbf 1_{\{X_u\neq v\}}b(X_u)\,du+\int_0^t
\mathbf 1_{\{X_u\neq v\}}
\sigma(X_u)\,dW_u^X
+L_t^0(R^X).\]
The local integrability conditions follow
from the corresponding assumptions for
$Y$ and the same change-of-variables
formula, after radial localization. Furthermore,
\[
\int_0^t\mathbf 1_{\{X_u=v\}}\,du
=\int_0^{S_t}
\mathbf 1_{\{Y_u=v\}}\,dC_u\\
=0.
\]
Thus $X$ spends zero Lebesgue time
at the vertex.

By Lemma~\ref{lem:excursion-preservation},
the excursions of $X$ and $Y$ away
from the vertex are in one-to-one
correspondence, with the same angular
marks. Moreover, for every
$\Lambda\in\mathcal B([0,2\pi))$,
the process
\(R_\Lambda^X(t)=R_\Lambda^Y(S_t)\)
is a continuous semimartingale,
by the standard time-change theorem
for semimartingales, applied locally
before explosion. Furthermore, by
Lemma~\ref{lem:local-time-partition-preservation},
for every $\Lambda\in\mathcal B([0,2\pi))$,
\[
L_t^0(R_\Lambda^X)
=
\nu(\Lambda)L_t^0(R^X),
\qquad 0\le t<\zeta_C,
\]
almost surely.

Finally, we identify the lifetime of $X$. For each $n$ such that the radial level $n$
is reached before $\tau^Y$, the strict
monotonicity of $C$ and the identity
$R_u^X=R_{S_u}^Y$ imply
\[
\tau_n^X=C_{\tau_n^Y}.
\]
If $\tau^Y<\infty$, radial explosion gives
$\tau_n^Y\uparrow\tau^Y$. Therefore,
\[
\tau^X
=\lim_{n\to\infty}C_{\tau_n^Y}
=\lim_{t\uparrow\tau^Y}C_t
=\zeta_C.
\]

Suppose next that $\tau^Y=\infty$.
If $\zeta_C=\infty$, then $X$ is defined
for all times and is nonexplosive.
It remains to consider the case
$\zeta_C<\infty$.

We first show that the radial path of $X$
cannot remain bounded on $[0,\zeta_C)$. For $n,m\in\mathbb N$, define
\[
E_{n,m}
:=
\left\{
\tau^Y=\infty,\quad
\zeta_C\leq m,\quad
\sup_{0\leq u<\zeta_C}R_u^X\leq n
\right\}.
\]
We claim that
\[
\mathbb P(E_{n,m})=0.
\]
Indeed, the sticky occupation identity
and the preservation of local time yield
\[
S_u=u+\eta L_u^0(R^X),
\qquad 0\leq u<\zeta_C.
\]
Since $C$ is continuous and strictly
increasing on $[0,\infty)$ and
$C_t\uparrow\zeta_C<\infty$ as $t\to\infty$,
we have
\[
S_u\longrightarrow\infty
\qquad\text{as }u\uparrow\zeta_C.
\]
Consequently, on $E_{n,m}$,
\[
L_u^0(R^X)
=\frac{S_u-u}{\eta}
\longrightarrow\infty
\qquad\text{as }u\uparrow\zeta_C.
\]
By assumption \textup{(LB)}, we can use
It\^o's formula for $e^{R^X}$ to obtain,
locally before $\zeta_C$,
\[
e^{R_u^X}
=e^{R_0^X}+N_u+L_u^0(R^X)+\int_0^u\mathbf1_{\{X_s\neq v\}}e^{R_s^X}
\left(b(X_s)+\frac12\sigma^2(X_s)\right)ds,
\]
where
\(N_u:=\int_0^u e^{R_s^X}\,dM_s^X\)
is a continuous local martingale. The coefficient of
the local-time term is $e^0=1$.
On $E_{n,m}$, $R_s^X\leq n$ and $u<\zeta_C\leq m$;
thus \textup{(LB)} gives
\[
N_u\leq e^n-e^{R_0^X}+K_ne^nm-L_u^0(R^X),
\qquad 0\leq u<\zeta_C.
\]
Since $L_u^0(R^X)\to\infty$ as $u\uparrow\zeta_C$,
we obtain
\[
N_u\longrightarrow-\infty
\qquad\text{as }u\uparrow\zeta_C.
\]
This is impossible for a continuous
local martingale on a stochastic interval
with finite endpoint.

Indeed, by the
Dambis--Dubins--Schwarz representation,
if its quadratic variation has a finite
limit at that endpoint, the martingale
converges to a finite limit. If its
quadratic variation diverges, Brownian
oscillation prevents convergence to
$-\infty$. Thus
\[
\mathbb P(E_{n,m})=0
\qquad\text{for all }n,m\in\mathbb N.
\]
Since
\[
\left\{
\tau^Y=\infty,\quad
\zeta_C<\infty,\quad
\sup_{u<\zeta_C}R_u^X<\infty
\right\} =
\bigcup_{n,m\geq1}E_{n,m},
\]
we conclude that, almost surely,
\[
\tau^Y=\infty,\quad \zeta_C<\infty
\quad\Longrightarrow\quad
\sup_{u<\zeta_C}R_u^X=\infty.
\]
It follows that every sufficiently large
radial level is reached before $\zeta_C$.
Consequently, $\tau_n^Y<\infty$ for all
sufficiently large $n$. Since $\tau^Y=\infty$
and $Y$ has continuous paths on every
compact time interval, we have
$\tau_n^Y\uparrow\infty$.
Therefore,
\[
\tau^X
=\lim_{n\to\infty}\tau_n^X
=\lim_{n\to\infty}C_{\tau_n^Y}
=\lim_{t\to\infty}C_t
=\zeta_C.
\]
Thus, in all cases,
\(\tau^X=\zeta_C.\) Consequently, $X$ satisfies, up to its
explosion time, the radial semimartingale
equation, the non-sticky occupation
property, the local-time partition
associated with $\nu$, and the required
continuity, adaptation, and local
integrability conditions.

Hence $X$ is a non-sticky Walsh diffusion
with coefficients $(b,\sigma)$ and
spinning measure $\nu$, in the framework
of Karatzas and Yan
~\cite{KaratzasYan2019}.

\medskip
\noindent\textit{Step 2. Identification of the clock.}

By the localized sticky occupation identity, for every \(n\in\mathbb N\),
\[
(t\wedge\tau_n^Y)-C_{t\wedge\tau_n^Y}
=
\eta
L_t^0\!\left(R^Y(\cdot\wedge\tau_n^Y)\right).
\]
Hence, on \(\{t<\tau_n^Y\}\),
\[
t-C_t=\eta L_t^0(R^Y).
\]
Since
\[
\{t<\tau^Y\}
=
\bigcup_{n\ge1}\{t<\tau_n^Y\},
\]
we obtain, up to the explosion time,
\[
t-C_t=\eta L_t^0(R^Y).
\]
Since $C_{S_t}=t$,
\[
t=S_t-\eta L_{S_t}^0(R^Y).
\]
Hence, by
Lemma~\ref{lem:local-time-partition-preservation},
\[
S_t=t+\eta L_t^0(R^X),
\qquad 0\le t<\zeta_C.
\]
This completes the proof of assertion~(i).

\medskip
\noindent\textit{Assertion (ii).}

Conversely, let $X$ be a non-sticky Walsh
diffusion associated with $(b,\sigma,\nu)$,
with explosion time $\tau^X$. Define
\[
A_t=t+\eta L_t^0(R^X),
\qquad 0\le t<\tau^X,
\]
and set
\(\zeta_A:=\lim_{t\uparrow\tau^X}A_t.\) For $0\le t<\zeta_A$, define
\[
T_t:=\inf\{u<\tau^X:A_u>t\},
\qquad
Y_t:=X_{T_t}.
\]
Since $A$ is continuous and strictly
increasing, $T$ is a continuous inverse
time change consisting of stopping times.

Let $(\mathcal F_t^X)$ denote the
underlying filtration of $X$, and set
\(\mathcal H_t:=\mathcal F_{T_t}^X,
\; 0\le t<\zeta_A.\) Consequently, $Y$ is continuous and
adapted to $(\mathcal H_t)$.
 Moreover, the time change preserves the excursions and their angular
marks. Hence the angular component of \(Y\) remains constant during
each excursion away from the vertex. 
By the change-of-variables identity
\eqref{eq:sticky-change-of-variables},
the required local integrability of $b$
and $\sigma^2$ is preserved after localization. Writing
\(R_t^Y=R_{T_t}^X\)
and time-changing the radial semimartingale
equation of $X$, we obtain
\[
R_t^Y
=
R_0^Y
+
\int_0^t
\mathbf 1_{\{Y_r\neq v\}}b(Y_r)\,dr
+
M_t^Y
+
L_{T_t}^0(R^X),
\]
where
\(M_t^Y:=M_{T_t}^X,\)
and $M^X$ denotes the continuous local
martingale part of the radial equation
of $X$.

Since $A_t\ge t$, its inverse satisfies
$T_t\le t$. Thus, by the standard
time-change theorem for continuous local
martingales, applied locally before
explosion, $M^Y$ is a continuous local
martingale with respect to $(\mathcal H_t)$. Moreover,
\[
\begin{aligned}
\langle M^Y\rangle_t=\langle M^X\rangle_{T_t}=\int_0^{T_t}
\mathbf 1_{\{X_u\neq v\}}
\sigma^2(X_u)\,du=\int_0^t
\mathbf 1_{\{Y_r\neq v\}}
\sigma^2(Y_r)\,dr.
\end{aligned}
\]
Consequently, possibly after enlarging
the probability space, there exists a
Brownian motion $W^Y$ such that
\[
M_t^Y
=
\int_0^t
\mathbf 1_{\{Y_r\neq v\}}
\sigma(Y_r)\,dW_r^Y.
\]
By the standard time-change property of
semimartingale local times,
\(
L_t^0(R^Y)=L_{T_t}^0(R^X).
\)
Therefore
\[
R_t^Y
=
R_0^Y
+
\int_0^t
\mathbf 1_{\{Y_r\neq v\}}b(Y_r)\,dr
+
\int_0^t
\mathbf 1_{\{Y_r\neq v\}}
\sigma(Y_r)\,dW_r^Y
+
L_t^0(R^Y).
\]
Furthermore, combining
\eqref{eq:sticky-vertex-occupation}
with the standard time-change identity
for local times recalled above yields
the sticky occupation identity
\eqref{eq:sticky-occupation}.

For every $\Lambda\in\mathcal B([0,2\pi))$,
the process $R_\Lambda^X$ is a continuous
semimartingale by
\cite[Theorem~2.12(ii)]{KaratzasYan2019}.
Since $T$ is a continuous time change
consisting of stopping times, the standard
time-change theorem implies that
$R_\Lambda^Y=R_\Lambda^X\circ T$ is also a continuous
semimartingale.

Moreover, since $X$ has spinning measure
$\nu$, the local-time partition property
and the standard time-change identity yield
\[
L_t^0(R_\Lambda^Y)=L_{T_t}^0(R_\Lambda^X)
=\nu(\Lambda)L_{T_t}^0(R^X)=\nu(\Lambda)L_t^0(R^Y),
\qquad 0\le t<\zeta_A.
\]

Finally, we identify the lifetime of $Y$.
For each $n$ such that $\tau_n^X<\tau^X$,
the continuity and strict monotonicity
of $A$, together with
$R^Y=R^X\circ T$, imply
\[
\tau_n^Y=A_{\tau_n^X}.
\]
Suppose first that $\tau^X<\infty$.
By the definition of radial explosion,
$\tau_n^X\uparrow\tau^X$. Consequently,
\[
\tau^Y
=\lim_{n\to\infty}\tau_n^Y
=\lim_{n\to\infty}A_{\tau_n^X}
=\lim_{u\uparrow\tau^X}A_u
=\zeta_A.
\]
Suppose next that $\tau^X=\infty$.
Since $A_u=u+\eta L_u^0(R^X)\ge u$,
we have $\zeta_A=\infty$ and
$T_t\le t$ for every finite $t$.
By continuity,
\[
\sup_{0\le s\le t}R_s^Y
=
\sup_{0\le u\le T_t}R_u^X
<\infty.
\]
Thus $Y$ cannot explode in finite time,
and therefore $\tau^Y=\infty=\zeta_A$.

We conclude that, in all cases,
$\tau^Y=\zeta_A$.
Consequently, $Y$ satisfies all the
conditions of
Definition~\ref{def:sticky-Walsh-diffusion}
and is therefore a sticky Walsh diffusion
associated with $(b,\sigma,\nu,\eta)$.

\end{proof}

\begin{assumption}
\label{ass:KY}

The coefficients \(b\) and \(\sigma\) satisfy the following
conditions:
\begin{enumerate}
    \item[(i)]
    \[
    \sigma_\theta(r)\neq0,
    \qquad
    r>0,\quad \theta\in[0,2\pi).
    \]

    \item[(ii)] For every fixed \(\theta\in[0,2\pi)\), the functions
    \[
    r\longmapsto
    \frac{b_\theta(r)}{\sigma_\theta^2(r)}
    \qquad\text{and}\qquad
    r\longmapsto
    \frac{1}{\sigma_\theta^2(r)}
    \]
    are locally integrable on \((0,\infty)\).

    \item[(iii)] There exists \(\delta>0\) such that
    \[
    \sup_{\substack{0<r\le\delta\\
                    \theta\in[0,2\pi)}}
    \frac{1+|b_\theta(r)|}
         {\sigma_\theta^2(r)}
    <\infty.
    \]
\end{enumerate}
\end{assumption}

The next result shows that under
Assumption~\ref{ass:KY} the lower-bound
condition \textup{(LB)} is unnecessary in the
de-stickifying direction. Its proof applies the
generalized Freidlin--Sheu formula
\cite[Theorem~2.12]{KaratzasYan2019}
directly to the sticky process $Y$.

\begin{prop}
\label{prop:desticky-KY}
Under Assumption~\ref{ass:KY}, the conclusion
of assertion~(i) of
Theorem~\ref{thm:characterization}
remains valid without assuming \textup{(LB)}.
\end{prop}

\begin{proof}
The proof of assertion~(i) of
Theorem~\ref{thm:characterization}, up to the
identification of the radial explosion time,
does not use \textup{(LB)}. The only case
requiring an alternative argument is
$\tau^Y=\infty$ and $\zeta_C<\infty$.

For each $\theta$, let
\[
s_\theta(r):=\int_0^r
\exp\left(-2\int_0^y
\frac{b_\theta(z)}{\sigma_\theta^2(z)}\,dz\right)dy.
\]
Under Assumption~\ref{ass:KY}, $s_\theta$ is
finite and increasing on compact radial intervals,
$s_\theta(0)=0$, and $s_\theta'(0+)=1$.
Moreover,
\[
b_\theta s_\theta'
+\frac12\sigma_\theta^2 s_\theta''=0
\qquad\text{a.e. on }(0,\infty).
\]
Fix $n\geq1$, and choose a smooth radial cutoff
$\chi_n$ equal to one on $[0,n+1]$ and zero
on $[n+2,\infty)$. Define
\[
F_n(r,\theta)
:=\chi_n(r)\frac{s_\theta(r)}{1+s_\theta(n+1)},
\qquad F_n(v)=0.
\]
The local integrability assumptions in \textup{(KY)}
and the uniform bound near the vertex in
Assumption~\ref{ass:KY}(iii) ensure that
$F_n$ belongs to the class of test functions
in \cite[Definition~2.6]{KaratzasYan2019}:
its raywise first derivative is locally
absolutely continuous, its vertex derivatives
are bounded, and its second derivatives admit
a common integrable bound near the vertex.
Also, for $0\leq r\leq n+1$,
$0\leq F_n(r,\theta)\leq1$.
Set
\[
c_n:=\int_0^{2\pi}
\frac{\nu(d\theta)}{1+s_\theta(n+1)}>0.
\]

By Definition~\ref{def:sticky-Walsh-diffusion},
$Y$ has continuous paths in the tree topology,
its radial part is the Skorokhod reflection of
the semimartingale driver appearing in its radial
equation, and its local times satisfy the
partition property. Thus, after radial
localization, it meets the hypotheses of the
generalized Freidlin--Sheu formula
\cite[Theorem~2.12]{KaratzasYan2019}.
That result applies to semimartingales on rays
without requiring non-stickiness.
Consequently, the scale equation yields,
up to the radial exit time $\tau_{n+1}^Y$
and locally before $\tau^Y$,
\[
F_n(Y_{t\wedge\tau_{n+1}^Y})
=F_n(Y_0)+N_{t\wedge\tau_{n+1}^Y}^{(n)}
+c_n L_{t\wedge\tau_{n+1}^Y}^0(R^Y),
\]
where $N^{(n)}$ is a continuous local martingale.

\noindent{Consider} the event
\[
E_n^Y:=\left\{\tau^Y=\infty,\ \zeta_C<\infty,
\ \sup_{t\geq0}R_t^Y\leq n\right\}.
\]
On this event, $\tau_{n+1}^Y=\infty$
and $F_n(Y_t)\leq1$ for all $t$.
The sticky occupation identity gives
\[
\eta L_t^0(R^Y)=t-C_t\longrightarrow\infty
\qquad(t\to\infty),
\]
because $C_t\uparrow\zeta_C<\infty$.
It follows that
\[
N_t^{(n)}\leq1-F_n(Y_0)-c_nL_t^0(R^Y)
\longrightarrow-\infty,
\]
which is impossible for a continuous local
martingale by the Dambis--Dubins--Schwarz
argument used in the proof of the theorem.
Hence $\mathbb P(E_n^Y)=0$ for every $n$.

Since $C$ is continuous and strictly increasing
and $X_t=Y_{S_t}$, the radial ranges of
$X$ and $Y$ agree before their respective
terminal times. Therefore, almost surely on
$\{\tau^Y=\infty,\zeta_C<\infty\}$,
the radial path of $X$ cannot remain bounded.
The correspondence of radial hitting times
and the continuity of $Y$ on finite intervals
then imply $\tau^X=\zeta_C$, exactly as
in the final part of Step~1 of the theorem.
\end{proof}

Assumption~\ref{ass:KY} is precisely
Condition~3.11 of Karatzas and
Yan~\cite{KaratzasYan2019}, specialized
to the full Walsh state space
$\mathcal W$, corresponding to
$\ell(\theta)\equiv\infty$.

Consequently, the existence and uniqueness result
\cite[Theorem~3.14]{KaratzasYan2019}, together
with Theorem~\ref{thm:characterization}(ii) and
Proposition~\ref{prop:desticky-KY}, yields the
following corollary.

\begin{corollary}
\label{cor:existence-uniqueness}

Suppose that Assumption~\ref{ass:KY} holds. Then, for every probability measure $\nu$ on
$[0,2\pi)$, every $\eta>0$, and every initial
probability distribution $\mu$ on $\mathcal W$,
there exists a sticky Walsh diffusion associated
with $(b,\sigma,\nu,\eta)$ and initial
distribution $\mu$. Moreover, its law is uniquely determined up to
the corresponding explosion time.

\end{corollary}

\begin{proof}
By \cite[Theorem~3.14]{KaratzasYan2019},
there exists a non-sticky Walsh diffusion
$X$ associated with $(b,\sigma,\nu)$ and
initial distribution $\mu$, unique in
distribution up to its explosion time $\tau^X$.
Define the clock $A$ as in \eqref{eq:sticky-clock}, with $R=\|X\|$, and let
\[
\zeta_A=\lim_{t\uparrow\tau^X}A_t,
\qquad
T_t=\inf\{s<\tau^X:A_s>t\},
\quad 0\leq t<\zeta_A.
\]
By Theorem~\ref{thm:characterization},
the process
\[
Y_t=X_{T_t},
\qquad 0\leq t<\zeta_A,
\]
is a sticky Walsh diffusion associated with
$(b,\sigma,\nu,\eta)$, with initial distribution
$\mu$ and explosion time $\tau^Y=\zeta_A$.
This proves existence.

Conversely, let $Y$ be any sticky Walsh
diffusion associated with $(b,\sigma,\nu,\eta)$
and initial distribution $\mu$. By
Proposition~\ref{prop:desticky-KY},
its de-stickifying time change produces a
non-sticky Walsh diffusion $X$ associated
with $(b,\sigma,\nu)$ and initial distribution
$\mu$.

The law of $X$ is uniquely determined by
\cite[Theorem~3.14]{KaratzasYan2019}. Since the radial local time admits a measurable
path-functional version, the law of
\[
\bigl(X,L^0(\|X\|)\bigr)
\]
is uniquely determined by the law of $X$.
Consequently, the additive clock $A$ defined in
\eqref{eq:sticky-clock} and its inverse determine the law of the
time-changed process $Y$ uniquely. This proves uniqueness in distribution up to
the corresponding explosion time.
\end{proof}

\section{The Wentzell operator and martingale problem}
\label{sec:wentzell-martingale}
The characterization theorem established in the previous section
provides a probabilistic construction of sticky Walsh diffusions
through a time-change correspondence with non-sticky Walsh diffusions. The purpose of this section is
to identify the differential operator naturally associated with this
construction and to derive the corresponding martingale problem.

Unlike the classical generator approach, the operator introduced below
is not used to define the process. Instead, it is derived from the
defining properties of sticky Walsh diffusions through the generalized
Freidlin--Sheu formula. This point of view allows us to derive a local Dynkin formula
and provides the analytical tools needed to study the local
behavior of sticky Walsh diffusions near the vertex. Under additional global regularity assumptions, we also
characterize the full generator domain on a suitable Banach space
and obtain an explicit representation of the resolvent.

\subsection{The Wentzell operator and test functions}

For every ray $\theta\in[0,2\pi)$,
define the second-order differential
operator by
\[
\mathcal{L}f(r,\theta)
:=\mathcal{L}_\theta f_\theta(r)
=
b_\theta(r)\partial_r f_\theta(r)
+
\frac12\sigma_\theta^2(r)\partial_{rr}f_\theta(r),
\qquad r>0.
\]
We now introduce a function space adapted
to the radial geometry of the Walsh state space
$\mathcal W$.

\begin{definition}[Uniform radial vanishing]
\label{def:radial-vanishing}
We denote by
\(\mathcal C_{\mathrm{rad},0}(\mathcal W)\)
the class of bounded functions
\(f:\mathcal W\to\mathbb R\) satisfying the
following conditions:
\begin{enumerate}
\item The function \(f\) is continuous with
respect to the tree topology.

\item The mapping
\((r,\theta)\longmapsto f_\theta(r)\)
is jointly Borel measurable on
\([0,\infty)\times[0,2\pi)\).

\item The function \(f\) vanishes uniformly
at radial infinity, namely,
\[
\lim_{R\to\infty}
\sup_{\theta\in[0,2\pi)}
\sup_{r\ge R}|f_\theta(r)|=0.
\]
\end{enumerate}
\end{definition}

\begin{remark}
The space
\(\mathcal C_{\mathrm{rad},0}(\mathcal W)\)
should not be confused with the usual
topological space \(C_0(\mathcal W)\),
defined through compact superlevel sets.
Since the Walsh state space has uncountably many
rays and is not locally compact at its
vertex, these two notions are different.
\end{remark}

\begin{definition}[Wentzell test functions]
\label{def:wentzell-domain}

We denote by \(\mathcal D_W(\mathcal{L})\) the class
of functions
\(f\in\mathcal C_{\mathrm{rad},0}(\mathcal W)\)
satisfying the following properties:

\begin{enumerate}

\item For every \(\theta\in[0,2\pi)\), the radial function
\(f_\theta(r):=f(r,\theta)\)
belongs to \(C^2((0,\infty))\cap C([0,\infty))\), and the one-sided radial derivative
\(\partial_r f_\theta(0+)\)
exists.

\item The first and second radial derivatives are locally bounded
uniformly in the angular variable: for every \(R>0\),
\[
\sup_{\theta\in[0,2\pi)}
\sup_{0<r\le R}
\left(
|\partial_r f_\theta(r)|
+
|\partial_{rr}f_\theta(r)|
\right)
<\infty.
\]

\item
The function $\mathcal{L}f$, defined on
$\mathcal W\setminus\{v\}$ as above,
admits a common limit $\mathcal{L}f(v)$
at the vertex, uniformly in the angular
variable, namely,
\[
\lim_{r\downarrow0}
\sup_{\theta\in[0,2\pi)}
|\mathcal{L}f(r,\theta)-\mathcal{L}f(v)|=0.
\]
Moreover, the resulting extension satisfies
$\mathcal{L}f\in\mathcal C_{\mathrm{rad},0}(\mathcal W)$.

\item The Wentzell boundary condition
\begin{equation}\label{eq:wentzell-condition}
\eta \mathcal{L}f(v)
=
\int_0^{2\pi}
\partial_r f_\theta(0+)\,\nu(d\theta)
\end{equation}
is satisfied.

\end{enumerate}

\end{definition}

\begin{remark}

The class $\mathcal D_W(\mathcal{L})$ is first introduced as a class of
test functions for the martingale problem under our general assumptions.
We do not identify it with the full generator domain under those
assumptions. Under the stronger global hypotheses of
Subsection~\ref{subsec:complete-generator}, we characterize the full
generator domain on a more restrictive Banach space.

\end{remark}

\subsection{A local Wentzell boundary value problem}

We next establish a local boundary value problem that will be used
in the analysis of the process near the vertex.

The following proposition extends Proposition~A.1 of Berry and
Colantoni \cite{BerryColantoni2026} from finite star graphs to the
Walsh state space with an arbitrary spinning probability measure.
In this setting, the boundary condition involves an integral over
the angular variable, and the construction requires joint
measurability and uniform estimates with respect to the direction.

\begin{prop}
\label{prop:local-Wentzell-BVP}

Let $\delta>0$, and consider the
raywise differential operators
$\mathcal L_\theta$ restricted to
$(0,\delta)$. Assume that, for every
$\theta\in[0,2\pi)$, the functions
$r\mapsto b_\theta(r)$ and
$r\mapsto\sigma_\theta(r)$ are continuous
on $[0,\delta]$, and that
\[
\sup_{\theta\in[0,2\pi)}
\sup_{0\le r\le\delta}
\left(
|b_\theta(r)|
+
|\sigma_\theta(r)|^{-1}
\right)<\infty.
\]
Let $f=(f_\theta)_{\theta\in[0,2\pi)}$
be bounded, jointly Borel measurable in
$(r,\theta)$, and continuous on $\mathcal W$
with respect to the tree topology.

Then there exists a unique family
$u=(u_\theta)_{\theta\in[0,2\pi)}$,
jointly Borel measurable in $(r,\theta)$,
such that
\(u_\theta\in
C^2((0,\delta))\cap C^1([0,\delta])\)
for every $\theta\in[0,2\pi)$, with
uniformly bounded first and second
radial derivatives, satisfying
\[
\mathcal{L}_\theta u_\theta=f_\theta,
\qquad 0<r<\delta,
\]
\[
u_\theta(\delta)=0,
\qquad
u_\theta(0)=u_{\theta'}(0),
\quad
\theta,\theta'\in[0,2\pi),
\]
and the Wentzell boundary condition
\[
\eta f(v)
=
\int_0^{2\pi}
u_\theta'(0+)\,\nu(d\theta).
\]
Moreover, suppose that the preceding
continuity and uniform boundedness
assumptions hold on a fixed neighborhood
$[0,\delta_0]$ of the vertex.
Then, as $\delta\downarrow0$,
\[
\sup_{\theta\in[0,2\pi)}
\sup_{0<r<\delta}
|u_\theta'(r)|=O(1),
\]
and
\[
\sup_{\theta\in[0,2\pi)}
\sup_{0<r<\delta}
|u_\theta''(r)|=O(1).
\]

\end{prop}

\begin{proof}

Set
\[
q_\theta(r)
:=
\frac{2b_\theta(r)}{\sigma_\theta^2(r)},
\qquad
h_\theta(r)
:=
\frac{2f_\theta(r)}{\sigma_\theta^2(r)},
\]
and
\[
Q_\theta(r)
:=
\int_0^r q_\theta(s)\,ds.
\]
Writing $p_\theta=u_\theta'$, the equation
$\mathcal{L}_\theta u_\theta=f_\theta$ becomes
\[
p_\theta'(r)+q_\theta(r)p_\theta(r)
=
h_\theta(r).
\]
By the integrating-factor method,
\[
p_\theta(r)
=
e^{-Q_\theta(r)}
\left(
a_\theta+
\int_0^r
h_\theta(z)e^{Q_\theta(z)}\,dz
\right),
\]
where
\(a_\theta=p_\theta(0)=u_\theta'(0+).\) Since $u_\theta(\delta)=0$, integration gives
\[
u_\theta(r)
=
-\int_r^\delta p_\theta(y)\,dy.
\]
Define
\[
\alpha_\theta(r)
:=
\int_r^\delta e^{-Q_\theta(y)}\,dy\quad\mbox{and}\quad
\beta_\theta(r)
:=
\int_r^\delta
e^{-Q_\theta(y)}
\left[
\int_0^y
h_\theta(z)e^{Q_\theta(z)}\,dz
\right]dy.
\]
Then
\[
u_\theta(r)
=
-a_\theta\alpha_\theta(r)
-\beta_\theta(r).
\]
Let $c$ denote the common value of $u$
at the vertex. The condition
$u_\theta(0)=c$ yields
\[
c
=
-a_\theta\alpha_\theta(0)
-\beta_\theta(0).
\]
Since $\alpha_\theta(0)>0$, we obtain
\[
a_\theta
=
-\frac{c+\beta_\theta(0)}
{\alpha_\theta(0)}.
\]
The Wentzell boundary condition becomes
\[
\eta f(v)
=
-\int_0^{2\pi}
\frac{c+\beta_\theta(0)}
{\alpha_\theta(0)}
\,\nu(d\theta).
\]
Consequently,
\[
c
=
-
\frac{
\displaystyle
\eta f(v)
+
\int_0^{2\pi}
\frac{\beta_\theta(0)}
{\alpha_\theta(0)}
\,\nu(d\theta)
}{
\displaystyle
\int_0^{2\pi}
\frac{1}{\alpha_\theta(0)}
\,\nu(d\theta)
}.
\]
We first verify that these expressions
are well defined. By the uniform boundedness assumptions,
there exists a constant $K_\delta<\infty$
such that
\[
|q_\theta(r)|\le K_\delta,
\qquad
|h_\theta(r)|\le K_\delta
\]
for all $\theta\in[0,2\pi)$ and
$0\le r\le\delta$. Hence \(
|Q_\theta(r)|\le K_\delta\delta,
\) and therefore
\[
\delta e^{-K_\delta\delta}
\le
\alpha_\theta(0)
\le
\delta e^{K_\delta\delta}.
\]
In particular, $\alpha_\theta(0)$
is strictly positive and bounded away
from zero uniformly in $\theta$. Moreover,
\[
|\beta_\theta(0)|
\le
\frac12 K_\delta
e^{2K_\delta\delta}\delta^2.
\]
Thus the integrals appearing in the
formula for $c$ are finite, and
\[
\int_0^{2\pi}
\frac{1}{\alpha_\theta(0)}
\,\nu(d\theta)>0.
\]
It follows that $c$ and every
$a_\theta$ are uniquely determined.

The preceding formulas therefore
define a family $u$ satisfying the
differential equation and all the
prescribed boundary conditions.

Since the coefficients and the forcing
term are jointly Borel measurable,
the functions
\((\theta,r)\longmapsto
Q_\theta(r),\;
\alpha_\theta(r),\;
\beta_\theta(r)\)
are jointly Borel measurable. Consequently, $\theta\mapsto a_\theta$
is Borel measurable, and the explicit
representation of $u_\theta(r)$ shows
that $(r,\theta)\mapsto u_\theta(r)$
is jointly Borel measurable. 

Furthermore, the continuity of
$b_\theta$, $\sigma_\theta$, and
$f_\theta$ along each ray implies
that $q_\theta$ and $h_\theta$
are continuous on $[0,\delta]$.
The representation of $p_\theta$
then gives
\(u_\theta\in
C^2((0,\delta))\cap C^1([0,\delta]).\)

We now establish the uniform bounds. For fixed $\delta$, the preceding
estimates imply that
$\alpha_\theta(0)^{-1}$ and
$\beta_\theta(0)$ are uniformly bounded.
Consequently, $c$ is finite and
\(
\sup_{\theta\in[0,2\pi)}|a_\theta|<\infty.
\) The representation
\[
u_\theta'(r)
=
e^{-Q_\theta(r)}
\left(
a_\theta+
\int_0^r
h_\theta(z)e^{Q_\theta(z)}\,dz
\right)
\]
therefore yields
\[
\sup_{\theta\in[0,2\pi)}
\sup_{0<r<\delta}
|u_\theta'(r)|<\infty.
\]
Since
\[
u_\theta''(r)
=
\frac{2}{\sigma_\theta^2(r)}
\left(
f_\theta(r)
-
b_\theta(r)u_\theta'(r)
\right),
\]
we also obtain
\[
\sup_{\theta\in[0,2\pi)}
\sup_{0<r<\delta}
|u_\theta''(r)|<\infty.
\]
The uniform bound on the first
derivatives implies continuity at
the vertex. Indeed,
\[
|u_\theta(r)-c|
\le
r
\sup_{\varphi\in[0,2\pi)}
\sup_{0<s<\delta}
|u_\varphi'(s)|.
\]
Hence
\[
\lim_{r\downarrow0}
\sup_{\theta\in[0,2\pi)}
|u_\theta(r)-c|=0.
\]
Thus $u$ is continuous at the vertex
with respect to the tree topology. Moreover, since
\(\mathcal{L}_\theta u_\theta(r)=f_\theta(r)\)
and $f$ is continuous at the vertex,
$\mathcal{L}u$ admits the common vertex value
\(\mathcal{L}u(v)=f(v),\)
with convergence uniform in $\theta$.

To prove uniqueness, let $u$ and
$\widehat u$ be two solutions and set
$w=u-\widehat u$. Then
\[
\mathcal{L}_\theta w_\theta=0,
\qquad
w_\theta(\delta)=0,
\]
and $w$ has a common vertex value $d$,
with
\[
\int_0^{2\pi}
w_\theta'(0+)\,\nu(d\theta)=0.
\]
The preceding representation, with
zero forcing term, gives
\[
w_\theta'(0+)
=
-\frac{d}{\alpha_\theta(0)}.
\]
Therefore,
\[
0
=
-d
\int_0^{2\pi}
\frac{1}{\alpha_\theta(0)}
\,\nu(d\theta).
\]
Since the integral is strictly positive,
we obtain $d=0$, and consequently
$w_\theta'(0+)=0$ for every $\theta$. The differential equation then yields
$w_\theta\equiv0$ on $[0,\delta]$.
This proves uniqueness.

Finally, assume that the coefficients
and $f$ satisfy the uniform bounds
on the fixed neighborhood
$[0,\delta_0]$. 
Since $f$ is bounded, the bounds on
$q_\theta$ and $h_\theta$ can be chosen
independently of $0<\delta\le\delta_0$. In particular, as $\delta\downarrow0$,
\[
\alpha_\theta(0)\asymp\delta
\]
uniformly in $\theta$, and
\[
\sup_\theta|\beta_\theta(0)|
=O(\delta^2).
\]
The explicit formula for $c$ gives
\(c=O(\delta).\) Consequently,
\(\sup_\theta|a_\theta|=O(1).\) Since
\[
\sup_{\theta}
\sup_{0\le r\le\delta}
|Q_\theta(r)|=O(\delta),
\]
the representation of $u_\theta'$
implies
\[
\sup_{\theta}
\sup_{0<r<\delta}
|u_\theta'(r)|=O(1).
\]
Finally, using
\[
u_\theta''(r)
=
\frac{2}{\sigma_\theta^2(r)}
\left(
f_\theta(r)
-
b_\theta(r)u_\theta'(r)
\right),
\]
we conclude that
\[
\sup_{\theta}
\sup_{0<r<\delta}
|u_\theta''(r)|=O(1).
\]
This completes the proof.
\end{proof}

\subsection{The Freidlin--Sheu formula and the martingale problem}
\label{subsec:martingale-problem}

We first recall the generalized Freidlin--Sheu
formula of Karatzas and Yan
\cite[Theorem~2.12]{KaratzasYan2019}
and verify the regularity conditions needed for its
application to the martingale problem for sticky Walsh diffusions.

Let $X=(R,\Theta)$ be a non-sticky Walsh
diffusion associated with $(b,\sigma,\nu)$,
defined up to its explosion time
\[
\tau^X=\lim_{n\to\infty}\tau_n^X,
\qquad
\tau_n^X
=
\inf\{t\ge0:R_t\ge n\}.
\]
Let $f\in\mathcal D_W(\mathcal{L})$.
For each $n\ge1$, 
choose a smooth radial cutoff $\chi_n$ such that
\[
\chi_n(r)=\begin{cases}
1, & 0\le r\le n,\\
0, & r\ge n+1,
\end{cases}
\]
and define
\(f^{(n)}(r,\theta)
=
\chi_n(r)f_\theta(r).\) By the definition of $\mathcal D_W(\mathcal{L})$,
the radial derivatives of $f^{(n)}$
are locally absolutely continuous
on $(0,\infty)$, and the first radial
derivatives have compact radial support. Moreover,
\[
\sup_{\theta\in[0,2\pi)}
|\partial_r f_\theta(0+)|<\infty,
\]
and, for every $a>0$,
\[
\sup_{\theta\in[0,2\pi)}
\sup_{0<r\le a}
|\partial_{rr}f_\theta(r)|<\infty.
\]
Thus the cutoff functions have the
radial regularity and boundedness
properties required for the generalized
Freidlin--Sheu formula
\cite[Definition~2.6,
Theorem~2.12]{KaratzasYan2019}. Applying that formula to $f^{(n)}$
and stopping at $\tau_n^X$, we obtain
\[f(X_{t\wedge\tau_n^X})-f(X_0)
=M_{t\wedge\tau_n^X}^{X,f}+\int_0^{t\wedge\tau_n^X}
\mathbf1_{\{R_s>0\}}\mathcal{L}f(X_s)\,ds+
\left(
\int_0^{2\pi}
\partial_r f_\theta(0+)\,\nu(d\theta)
\right)
L_{t\wedge\tau_n^X}^0(R),
\]
where $M^{X,f}$ is a continuous local
martingale on the stochastic interval
$[0,\tau^X)$, given by
\[
M_t^{X,f}
=
\int_0^t
\mathbf1_{\{R_s>0\}}
\sigma(X_s)
\partial_r f_{\Theta_s}(R_s)\,dW_s,
\qquad 0\le t<\tau^X.
\]
Here $L^0(R)$ denotes the radial local
time with the normalization used in
the Skorokhod decomposition
\[
R_{t\wedge\tau_n^X}
=R_0
+\int_0^{t\wedge\tau_n^X}
\mathbf1_{\{R_s>0\}}b(X_s)\,ds+\int_0^{t\wedge\tau_n^X}
\mathbf1_{\{R_s>0\}}\sigma(X_s)\,dW_s+L_{t\wedge\tau_n^X}^0(R).
\]

\begin{prop}
\label{prop:sticky-martingale-problem}

Let $Y$ be a sticky Walsh diffusion
associated with $(b,\sigma,\nu,\eta)$,
with explosion time $\tau^Y$,
and let $f\in\mathcal D_W(\mathcal{L})$. Let $(\tau_n^Y)_{n\ge1}$ be the
radial exit times introduced in
Definition~\ref{def:sticky-Walsh-diffusion}. Then
\[
M_{t\wedge\tau_n^Y}^{Y,f}
:=
f(Y_{t\wedge\tau_n^Y})-f(Y_0)-\int_0^{t\wedge\tau_n^Y}
\mathcal{L}f(Y_s)\,ds
\]
is a continuous martingale. Consequently,
\[
M_t^{Y,f}
=
f(Y_t)-f(Y_0)
-\int_0^t\mathcal{L}f(Y_s)\,ds,
\qquad 0\le t<\tau^Y,
\]
is a continuous local martingale.

\end{prop}

\begin{proof}

Fix $n\geq1$ and set
$s_n(t):=t\wedge\tau_n^Y$.
The radial semimartingale equation and the
tree-topology continuity in
Definition~\ref{def:sticky-Walsh-diffusion}
show that, locally before explosion, $Y$
is a semimartingale on rays in the sense of
\cite[Definition~2.4]{KaratzasYan2019}.
The partition of local time in the same
definition permits the application of the
generalized Freidlin--Sheu formula
\cite[Theorem~2.12]{KaratzasYan2019}
directly to $Y$. As recalled above, the
required test-function regularity follows
by radial cutoff and localization.
Thus
\[
f(Y_{s_n(t)})-f(Y_0)
=N_{s_n(t)}^{Y,f}
+\int_0^{s_n(t)}\mathbf1_{\{Y_u\neq v\}}
\mathcal Lf(Y_u)\,du+\left(\int_0^{2\pi}
\partial_r f_\theta(0+)\,\nu(d\theta)\right)
L_{s_n(t)}^0(R^Y),
\]
where $N^{Y,f}$ is a continuous local martingale.
By the Wentzell boundary condition
\eqref{eq:wentzell-condition}, the coefficient
of the local-time term is $\eta\mathcal Lf(v)$.
The sticky occupation identity gives
\[
\eta L_{s_n(t)}^0(R^Y)
=\int_0^{s_n(t)}\mathbf1_{\{Y_u=v\}}\,du.
\]
Combining the finite-variation terms yields
\[
f(Y_{s_n(t)})-f(Y_0)
=N_{s_n(t)}^{Y,f}
+\int_0^{s_n(t)}\mathcal Lf(Y_u)\,du.
\]
Therefore
\[
M_{t\wedge\tau_n^Y}^{Y,f}=N_{t\wedge\tau_n^Y}^{Y,f}
\]
is a continuous local martingale.
Since $f,\mathcal Lf\in
\mathcal C_{\mathrm{rad},0}(\mathcal W)$,
both functions are bounded. Hence, on every
bounded time interval,
\[
|M_{t\wedge\tau_n^Y}^{Y,f}|
\leq2\|f\|_\infty+t\|\mathcal Lf\|_\infty.
\]
It follows that the stopped process is a true
martingale. Finally, since $\tau_n^Y\uparrow\tau^Y$,
the sequence $(\tau_n^Y)$ localizes $M^{Y,f}$
on $[0,\tau^Y)$. This completes the proof.
\end{proof}

\subsection{Local Dynkin's formula}
\label{subsec:Dynkin-formula}

The martingale problem established in
Proposition~\ref{prop:sticky-martingale-problem}
yields Dynkin's formula for global test functions.
For the local analysis in Section~4, we establish
a corresponding formula for the process stopped
upon exiting a radial neighborhood of the vertex,
under assumptions imposed only on that neighborhood.

For $\delta>0$, let
\[
\tau_\delta
:=
\inf\{t\in[0,\tau^Y):\|Y_t\|=\delta\}.
\]

\begin{prop}
\label{prop:local-Dynkin}

Let $\delta>0$ and set
\(\mathcal W_\delta
:=
\{y\in\mathcal W:\|y\|\le\delta\}.\)
Fix $x\in\mathcal W_\delta$ and let $Y$ be a sticky Walsh diffusion
associated with $(b,\sigma,\nu,\eta)$, started from $x$,
with $\mathbb E_x$ denoting expectation under its law.
Let $f$ be a function defined on $\mathcal W_\delta$.
Assume that $f$ is continuous with respect
to the tree topology and that
\((r,\theta)\longmapsto f_\theta(r)\)
is jointly Borel measurable on
$[0,\delta]\times[0,2\pi)$.
Suppose that, for every $\theta\in[0,2\pi)$,
\(f_\theta\in C^1([0,\delta])\cap C^2((0,\delta)),\)
and that
\[
\sup_{\theta\in[0,2\pi)}
\sup_{0<r<\delta}
\left(
|\partial_r f_\theta(r)|
+
|\partial_{rr}f_\theta(r)|
\right)<\infty.
\]
Assume moreover that $\mathcal Lf$ is bounded
on $\mathcal W_\delta$, admits a common
value $\mathcal Lf(v)$ at the vertex,
and satisfies the local Wentzell boundary
condition \eqref{eq:wentzell-condition}.
Then, for every $t\ge0$,
\begin{equation}\label{eq:stopped-Dynkin}
\mathbb E_x\bigl[f(Y_{t\wedge\tau_\delta})\bigr]
=
f(x)
+
\mathbb E_x\left[
\int_0^{t\wedge\tau_\delta}\mathcal Lf(Y_s)\,ds
\right].
\end{equation}
If, in addition, $\mathbb E_x[\tau_\delta]<\infty$,
then
\begin{equation}\label{eq:terminal-Dynkin}
\mathbb E_x[f(Y_{\tau_\delta})]
=
f(x)
+
\mathbb E_x\left[
\int_0^{\tau_\delta}\mathcal Lf(Y_s)\,ds
\right].
\end{equation}
\end{prop}

\begin{proof}

We first construct an extension of $f$
to the whole Walsh state space. By the uniform bound on the second
radial derivatives, there exists
$C<\infty$ such that
\[
\left|
\partial_r f_\theta(s)
-
\partial_r f_\theta(r)
\right|
\le C|s-r|,
\qquad 0<r<s<\delta,
\]
uniformly in $\theta$. Consequently, the one-sided limits
$\partial_r f_\theta(0+)$ and
$\partial_r f_\theta(\delta-)$ exist
and are finite. Moreover,
\[
\sup_{\theta\in[0,2\pi)}
\left|
\partial_r f_\theta(\delta-)
\right|<\infty.
\]
Choose $\chi\in C^\infty([0,\infty))$
such that
\[
\chi(s)=1
\quad\text{for }0\le s\le\frac13,
\qquad
\chi(s)=0
\quad\text{for }s\ge1.
\]
Define $\widetilde f$ by
\[
\widetilde f_\theta(r)
=
f_\theta(r),
\qquad 0\le r\le\delta,
\]
and, for $r\ge\delta$,
\[
\widetilde f_\theta(r)
=
\chi(r-\delta)
\left[
f_\theta(\delta)
+
(r-\delta)
\partial_r f_\theta(\delta-)
\right].
\]
Since the values and first radial
derivatives agree at $r=\delta$,
the extension is continuously
differentiable along each ray. Furthermore, its first radial
derivative is absolutely continuous
on every compact radial interval.
Indeed, it is absolutely continuous
on $[0,\delta]$ by the uniform bound
on the second derivatives, and it
is smooth on $(\delta,\infty)$, with
matching first derivatives at
$r=\delta$.

For each $r\in(0,\delta)$, the map
$\theta\mapsto\partial_r f_\theta(r)$
is Borel measurable, since the derivative
is the pointwise limit of measurable
difference quotients. Consequently,
$\theta\mapsto\partial_r f_\theta(\delta-)$
is Borel measurable, being the pointwise
limit of $\partial_r f_\theta(r_k)$
for any sequence $r_k\uparrow\delta$. It follows that
$(r,\theta)\mapsto\widetilde f_\theta(r)$
is jointly Borel measurable.

Since $f$ is bounded on
$\mathcal W_\delta$ and its first
radial derivatives are uniformly
bounded, the extension has bounded
first and second radial derivatives
away from the joining point,
uniformly in $\theta$.

In particular, the first radial derivative
of $\widetilde f_\theta$ is absolutely
continuous on $[0,\infty)$ for every
$\theta\in[0,2\pi)$, including across
the joining point $r=\delta$.
Furthermore, the first radial derivatives
at the vertex are uniformly bounded,
and the second radial derivatives admit
a uniform integrable bound near the vertex.

Moreover, $\widetilde f$ has compact
radial support. In particular,
$\widetilde f$ satisfies the
conditions of
\cite[Definition~2.6]{KaratzasYan2019}.

If $\|x\|=\delta$, then $\tau_\delta=0$
and the two identities are immediate.
We therefore assume $\|x\|<\delta$.
By continuity of the radial paths and the
definition of radial explosion,
$\tau_\delta\le\tau_n^Y$ for any integer
$n>\delta$. The generalized Freidlin--Sheu
formula \cite[Theorem~2.12]{KaratzasYan2019},
applied directly to the sticky process $Y$
and the extension $\widetilde f$, and stopped
at $t\wedge\tau_\delta$, gives
\[
f(Y_{t\wedge\tau_\delta})-f(x)=N_{t\wedge\tau_\delta}^{Y,f}
+\int_0^{t\wedge\tau_\delta}
\mathbf1_{\{Y_s\ne v\}}\mathcal Lf(Y_s)\,ds+\left(\int_0^{2\pi}
\partial_r f_\theta(0+)\,\nu(d\theta)\right)
L_{t\wedge\tau_\delta}^0(R^Y),
\]
where $N^{Y,f}$ is a continuous local martingale.
Indeed, the radial semimartingale equation
and local-time partition property in
Definition~\ref{def:sticky-Walsh-diffusion}
allow the application of this formula to $Y$,
after radial localization. The local Wentzell boundary condition and
the sticky occupation identity give
\[
\left(\int_0^{2\pi}
\partial_r f_\theta(0+)\,\nu(d\theta)\right)
L_{t\wedge\tau_\delta}^0(R^Y)
=
\int_0^{t\wedge\tau_\delta}
\mathbf1_{\{Y_s=v\}}\mathcal Lf(v)\,ds.
\]
Consequently,
\[
M_{t\wedge\tau_\delta}^{Y,f}
:=f(Y_{t\wedge\tau_\delta})-f(x)
-\int_0^{t\wedge\tau_\delta}\mathcal Lf(Y_s)\,ds
=N_{t\wedge\tau_\delta}^{Y,f}
\]
is a continuous local martingale.
Since $f$ and $\mathcal Lf$ are bounded
on $\mathcal W_\delta$, for every $T>0$
and $0\le t\le T$,
\[
|M_{t\wedge\tau_\delta}^{Y,f}|
\le 2\|f\|_{\infty,\delta}
+T\|\mathcal Lf\|_{\infty,\delta},
\]
where the norms are taken over $\mathcal W_\delta$.
Thus the stopped local martingale is a true
martingale on each bounded time interval.
Taking expectations yields
\eqref{eq:stopped-Dynkin}.

If $\mathbb E_x[\tau_\delta]<\infty$,
then $\tau_\delta<\infty$ almost surely.
Moreover, the boundedness of $f$ and
$\mathcal Lf$ on $\mathcal W_\delta$
implies
\[
\left|\int_0^{t\wedge\tau_\delta}
\mathcal Lf(Y_s)\,ds\right|
\le\tau_\delta\|\mathcal Lf\|_{\infty,\delta}.
\]
By dominated convergence,
$\mathbb E_x[f(Y_{t\wedge\tau_\delta})]$
converges to $\mathbb E_x[f(Y_{\tau_\delta})]$,
and the integral term in
\eqref{eq:stopped-Dynkin} converges in $L^1$
to its terminal-time counterpart.
Letting $t\to\infty$ yields
\eqref{eq:terminal-Dynkin}.
\end{proof}

\subsection{The complete generator under global regularity}
\label{subsec:complete-generator}

The Wentzell martingale problem above is formulated under the local
assumptions used in the probabilistic construction. Under stronger,
global assumptions on the coefficients, the same Wentzell condition
also characterizes the \emph{entire} infinitesimal-generator domain
on a natural space of uniformly tree-continuous functions.
This statement concerns the generator on a specified Banach space;
it does not identify a generator domain under the full standing
assumptions, nor on all bounded Borel functions.

\Needspace{10\baselineskip}
Throughout this subsection, we impose the following additional global regularity condition on the coefficients $b_\theta$ and $\sigma_\theta$.

\medskip
\noindent\textbf{Assumption (GR).}
There exist constants $B,C,K>0$ and $c>0$ such that, for all $r,s>0$ and $\theta\in S^1$,
\begin{align*}
& |b_\theta(r)|\le B,\qquad 0<c\le |\sigma_\theta(r)|\le C,\\
& |b_\theta(r)-b_\theta(s)|+|\sigma_\theta(r)-\sigma_\theta(s)|\le K|r-s|.
\end{align*}
The coefficients are understood to be extended continuously to $r=0$; these extensions remain jointly Borel. No continuity in $\theta$ is assumed.

\medskip
The sticky Walsh diffusion $Y$ is the process constructed in
Section~\ref{sec:probabilistic-construction}. Under
Assumption~\textup{(GR)} it is nonexplosive: the bounded radial
drift and diffusion coefficients prevent finite-time explosion of its
non-sticky counterpart, while the inverse sticky clock satisfies
$T_t\le t$. The generator argument below relies on the established
martingale problem and the discounted Dynkin identity
derived from it; no separate strong Markov assertion
for $Y$ is needed. Define
\[
\mathscr B:=
\left\{
f\in\mathcal C_{\mathrm{rad},0}(\mathcal W):
f\text{ is uniformly continuous with respect to }\rho
\right\},
\]
equipped with the supremum norm. This space is complete. 

For the resolvent estimates, let $\lambda>0$ and $g\in\mathscr B$.
For each $\theta\in S^1$, let $Z^\theta$ be the one-dimensional
diffusion with coefficients $(b_\theta,\sigma_\theta)$,
started from $r\ge0$, and let
$H_0:=\inf\{t\ge0:Z_t^\theta=0\}$ denote its first hitting time of zero.
The diffusion is killed at $H_0$. Define
\begin{align*}
h_{\lambda,\theta}(r)&=\mathbb E_r^\theta[e^{-\lambda H_0};H_0<\infty],\\
w_{\lambda,\theta}(r)&=\mathbb E_r^\theta\!\left[\int_0^{H_0}e^{-\lambda t}g_\theta(Z_t^\theta)\,dt\right].
\end{align*}
We first establish three auxiliary lemmas.

\begin{lem}
\label{gen:grad}
Let $\lambda>0$, and let $(g_\theta)_{\theta\in S^1}$
be a family of bounded functions on $[0,\infty)$.
Suppose that
$u_\theta\in C^2((0,\infty))\cap C_b([0,\infty))$
solves $(\lambda-\mathcal L_\theta)u_\theta=g_\theta$
and that
$\sup_\theta\|u_\theta\|_\infty\le H$,
$\sup_\theta\|g_\theta\|_\infty\le G$,
for some nonnegative constants $H$ and $G$.
Then
$\sup_{\theta\in S^1,\,r>0}
(|u_\theta'(r)|+|u_\theta''(r)|)<\infty$.
All $u_\theta'(0+)$ exist, and the
derivative bound is uniform in $\theta$.
\end{lem}

\begin{proof}
Write $p_\theta=u_\theta'$, $a_\theta=2b_\theta/\sigma_\theta^2$ and $q_\theta=2(\lambda u_\theta-g_\theta)/\sigma_\theta^2$. Then
\[
p_\theta'+a_\theta p_\theta=q_\theta
\]
on $(0,\infty)$ with
\[
|a_\theta|\le A:=2B/c^2,\qquad |q_\theta|\le D:=2(\lambda H+G)/c^2.
\]
For any $r>0$ and $0\le t\le1$, variation of constants gives
\[
p_\theta(r+t)=e^{-\int_0^t a_\theta(r+s)ds}p_\theta(r)+\int_0^t e^{-\int_s^t a_\theta(r+z)dz}q_\theta(r+s)ds.
\]
Integrate from $t=0$ to $1$. The coefficient of $p(r)$ is at least $e^{-A}$; the remaining term has absolute value at most $De^{A}$. As
\[
|u_\theta(r+1)-u_\theta(r)|\le2H,
\]
we obtain
\[
|p_\theta(r)|\le e^A(2H+De^A),
\]
for all $r,\theta$. The differential equation gives
\[
|u''_\theta|\le A|p|+D
\]
uniformly. Thus $u_\theta'$ is uniformly Lipschitz on $(0,\infty)$, and each $u_\theta'(0+)$ exists with the asserted bound.
\end{proof}

\begin{lem}\label{gen:dir}
The functions $h_{\lambda,\theta}$ and $w_{\lambda,\theta}$ are jointly Borel, bounded, and $C^2$ in $(0, \infty)$, with continuous boundary values $h(0)=1,w(0)=0$. They uniquely solve
$(\lambda-\mathcal{L}_\theta)h=0$ and $(\lambda-\mathcal{L}_\theta)w=g_\theta$ among bounded solutions with those boundary data. Their derivatives up to order two are bounded uniformly in $\theta,r$. There exist $0<\kappa_\lambda\le K_\lambda<\infty$ such that
\[
e^{-K_\lambda r}\le h_{\lambda,\theta}(r)\le e^{-\kappa_\lambda r},\qquad
-K_\lambda\le h_{\lambda,\theta}'(0+)\le-\kappa_\lambda.
\]
Moreover $\sup_\theta\sup_{r\ge R}|w_{\lambda,\theta}(r)|\to0$ as $R\to\infty$.
\end{lem}
\begin{proof}
Fix $\theta\in S^1$ and $\lambda>0$. Under Assumption~\textup{(GR)}, the diffusion $Z^\theta$ exists as a nonexplosive strong Markov process. Let $0<a<r<b<\infty$, and put
\[
\tau_{a,b}:=\inf\{t\ge0:Z_t^\theta\notin(a,b)\}.
\]
Uniform ellipticity on the bounded interval implies
$\mathbb P_r^\theta(\tau_{a,b}<\infty)=1$.
Since $H_0\ge\tau_{a,b}$ for $r\in(a,b)$, the strong Markov property yields
\[
 h_{\lambda,\theta}(r)
 =\mathbb E_r^\theta\!\left[
 e^{-\lambda\tau_{a,b}}
 h_{\lambda,\theta}(Z_{\tau_{a,b}}^\theta)\right]
\]
and
\[
 w_{\lambda,\theta}(r)=\mathbb E_r^\theta\!\left[
 \int_0^{\tau_{a,b}}e^{-\lambda t}
 g_\theta(Z_t^\theta)\,dt\right]+\mathbb E_r^\theta\!\left[
 e^{-\lambda\tau_{a,b}}
 w_{\lambda,\theta}(Z_{\tau_{a,b}}^\theta)\right].
\]
On $[a,b]$, $b_\theta$ and $\sigma_\theta$ are continuous,
$\sigma_\theta^2$ is bounded away from zero, and $g_\theta$ is
continuous because $g\in\mathscr B$.
The classical two-point Dirichlet problems with boundary values
$h_{\lambda,\theta}(a),h_{\lambda,\theta}(b)$ and
$w_{\lambda,\theta}(a),w_{\lambda,\theta}(b)$ therefore
have unique $C^2((a,b))\cap C([a,b])$ solutions. The
bounded-interval Feynman--Kac formula identifies these solutions
with the expectations above. As $a,b$ are arbitrary, it follows that
$h_{\lambda,\theta},w_{\lambda,\theta}\in C^2((0,\infty))$
and that
\[
(\lambda-\mathcal L_\theta)h_{\lambda,\theta}=0,
\qquad
(\lambda-\mathcal L_\theta)w_{\lambda,\theta}=g_\theta
\quad(r>0).
\]
Choose $\kappa_\lambda>0$ sufficiently small that
\[
C^2\kappa_\lambda^2/2+B\kappa_\lambda\le\lambda,
\]
and $K_\lambda$ sufficiently large that
\[
c^2 K_\lambda^2/2-BK_\lambda\ge\lambda.
\] The functions $e^{-\kappa_\lambda r}$ and $e^{-K_\lambda r}$ are respectively super- and subsolutions of $(\lambda-\mathcal{L}_\theta)h=0$. Stopped It\^o/optional sampling and passage to the limit give the exponential comparison. In particular, $h_{\lambda,\theta}(r)\to1$ as $r\downarrow0$, uniformly in $\theta$.

We also make the boundary convergence of $w$ uniform in the angular
variable. Since $g$ is bounded and $H_0$ is the killing time, the
Feynman--Kac representation and the exponential lower bound for $h$
give, for every $r\ge0$,
\[
\begin{aligned}
|w_{\lambda,\theta}(r)|
&\le \|g\|_\infty\,\mathbb E_r^\theta
     \left[\int_0^{H_0}e^{-\lambda t}\,dt\right]\\
&=\frac{\|g\|_\infty}{\lambda}
  (1-h_{\lambda,\theta}(r))\\
&\le\frac{\|g\|_\infty}{\lambda}
  (1-e^{-K_\lambda r}).
\end{aligned}
\]
Consequently,
$\sup_\theta|w_{\lambda,\theta}(r)|\to0$ as
$r\downarrow0$. This justifies the common boundary value $w(0)=0$
with the uniformity required for the tree topology.

For uniqueness among bounded solutions with these boundary values,
apply It\^o's formula to the difference of two solutions, stopped
at $H_\varepsilon\wedge H_R\wedge t$. Let $R\uparrow\infty$,
then $\varepsilon\downarrow0$, and finally $t\uparrow\infty$;
boundedness and $\lambda>0$ make the terminal term vanish.

For completeness, we spell out the parameter-measurability argument. Extend $b_\theta$ and $\sigma_\theta$ to the real line by setting their values at $x\le0$ equal to those at $0$; the extensions remain jointly Borel and uniformly Lipschitz in $x$. 
Let $Z^{\theta,r}$ denote the diffusion
$Z^\theta$ started from $r$.
On a common Wiener space, construct the
Euler--Maruyama approximations
$Z^{\theta,r,(n)}$ of the corresponding
one-dimensional SDE, using dyadic time steps. For each $n$, the map
\[
(\theta,r,\omega)\longmapsto
Z^{\theta,r,(n)}(\omega)|_{[0,T]}
\]
is Borel as a map into $C([0,T])$. The usual Lipschitz SDE estimates, with constants uniform in $\theta$ and in $r$ on bounded sets, give, for each $T,M<\infty$,
\[
\sup_{\theta\in S^1,\,0\le r\le M}
\mathbb E\!\left[
\sup_{t\le T}|Z_t^{\theta,r,(n)}-Z_t^{\theta,r}|^2
\right]\le C_{T,M}2^{-n}.
\]
Thus the Euler approximations converge almost surely uniformly on $[0,T]$ for each fixed $(r,\theta)$, and their Borel pointwise limit (defined arbitrarily on the null exceptional set) provides a jointly measurable version of the solution paths. Taking compatible integer horizons gives a version on $[0,\infty)$. For $r>0$ the first hitting time $H_0$ is a Borel functional of a continuous path. Since $g_\theta$ is jointly Borel and bounded, the maps
\[
(\theta,r,\omega)\longmapsto
 e^{-\lambda H_0}\mathbf1_{\{H_0<\infty\}},
\qquad
(\theta,r,\omega)\longmapsto
\int_0^{H_0}e^{-\lambda t}g_\theta(Z_t^{\theta,r})\,dt
\]
are jointly measurable; taking expectations proves joint Borel measurability of $h_{\lambda,\theta}(r)$ and $w_{\lambda,\theta}(r)$. The uniform bounds $\|h\|_\infty\le1$ and
$\|w\|_\infty\le\|g\|_\infty/\lambda$, together with the
boundary continuity just proved, allow us to apply
Lemma~\ref{gen:grad}. Thus the first two radial derivatives
are bounded uniformly in $\theta$ and $r$, and the derivatives
at $0+$ exist. From the exponential comparison,
\[
-K_\lambda\le h_{\lambda,\theta}'(0+)\le-\kappa_\lambda.
\]
The boundary derivatives are jointly measurable in $\theta$, since
\[
\partial_r h_{\lambda,\theta}(0+)
=\lim_{m\to\infty}m\bigl(h_{\lambda,\theta}(m^{-1})-1\bigr),\quad
\partial_r w_{\lambda,\theta}(0+)
=\lim_{m\to\infty}m\,w_{\lambda,\theta}(m^{-1}).
\]
Consequently, the angular integrals of these derivatives are well defined.

For the uniform tail of $w$, fix $R>0$ and $r\ge R$. Stop the killed diffusion at its first hitting time $H_R$. Before $H_R$ it lies in $[R,\infty)$, hence the integral contributes at most $\lambda^{-1}\sup_{\theta,s\ge R}|g_\theta(s)|$. The post-$H_R$ contribution is bounded by
\[
\|g\|_\infty\lambda^{-1}\mathbb E_r[e^{-\lambda H_R};H_R<\infty]\le \|g\|_\infty\lambda^{-1}e^{-\kappa_\lambda(r-R)}.
\] Given $\varepsilon>0$, choose $R_0$ so large that $\sup_{\theta,s\ge R_0}|g_\theta(s)|\le\varepsilon$; for $r\ge R_0+M$, the displayed bound becomes
\[
\varepsilon/\lambda+(\|g\|_\infty/\lambda)e^{-\kappa_\lambda M}.
\]
Now let $\varepsilon\downarrow0$ and $M\uparrow\infty$. This proves uniform radial vanishing.
\end{proof}

\begin{lem}\label{gen:usc}
There is a constant $C_0$, depending only on the global bounds on $b,\sigma$, such that for $t\ge0$,
\[
\sup_{x\in\mathcal W}\mathbb E_x\!\left[\sup_{s\le t}\rho(Y_s,x)\right]\le C_0(\sqrt t+t).
\]
\end{lem}
\begin{proof}
For the nonsticky process $X$, its radial component $R$ is the Skorokhod reflection of
\[
r+M_t+\int_0^t b(X_s)ds,
\]
where $\langle M\rangle_t\le C^2t$. The Skorokhod map and the Burkholder–Davis–Gundy inequality imply that, for
some constant $C_1>0$ independent of $x$ and $t$,
\[
\sup_{x\in\mathcal W}
\mathbb E_x\left[\sup_{s\le t}|R_s-r|\right]
\le C_1(\sqrt t+t).
\]
As long as the process remains on its initial ray,
the tree distance from $Y_s$ to $x$ equals
the radial distance. If the ray changes, the path has hit the vertex, which implies
\[
r\le\sup_{s\le t}|R_s-r|;
\]
hence
\[
\rho(X_s,x)=r+R_s\le3\sup_{u\le t}|R_u-r|.
\] Thus the same estimate, up to a factor 3, holds in tree distance. The inverse sticky clock satisfies $T_s\le s$; $Y_s=X_{T_s}$ implies the stated bound.
\end{proof}

\begin{thm}[Full generator domain on $\mathscr B$]\label{gen:gen}
Under Assumption~\textup{(GR)}, the transition operators $P_t f(x)=\mathbb E_x[f(Y_t)]$ form a strongly continuous contraction semigroup on $\mathscr B$. If $\mathcal L_{\mathscr B}$ denotes its infinitesimal generator on this Banach space, then
\[
\begin{aligned}
D(\mathcal L_{\mathscr B})=\biggl\{f\in\mathscr B:\;&
\exists\,F\in\mathscr B\text{ such that}\;f_\theta\in C^2((0,\infty)),\,\theta\in S^1,\\
&\mathcal L_\theta f_\theta(r)=F(r,\theta),
\; r>0,\ \theta\in S^1,\;\int_0^{2\pi}f_\theta'(0+)\,\nu(d\theta)=\eta F(v)
\biggr\}.
\end{aligned}
\]
For every $f$ in this domain,
$\mathcal L_{\mathscr B}f=F$. For any $\lambda>0$ and $g\in\mathscr B$, the full resolvent satisfies
\begin{equation}\label{eq:resolvent-representation}
\mathcal R_\lambda g(r,\theta)=w_{\lambda,\theta}(r)+h_{\lambda,\theta}(r)c_\lambda,
\end{equation}
where
\begin{equation}\label{eq:resolvent-vertex}
c_\lambda=\frac{\eta g(v)+\int_0^{2\pi}w_{\lambda,\theta}'(0+)\nu(d\theta)}{\eta\lambda-\int_0^{2\pi}h_{\lambda,\theta}'(0+)\nu(d\theta)}.
\end{equation}
\end{thm}

\begin{proof}[Proof of Theorem~\ref{gen:gen}]
\leavevmode\par\noindent
\textbf{Step 1 (constructing and identifying the resolvent).} For $g\in\mathscr B$, let $h,w$ be the functions defined above. Lemma~\ref{gen:dir} shows that the denominator $\eta\lambda-\int h'_{\lambda,\theta}(0+)\nu(d\theta)$ is at least $\eta\lambda+\kappa_\lambda>0$, and that all integrals are well defined. Define $c_\lambda$ by \eqref{eq:resolvent-vertex} and set
$u_\theta=w_{\lambda,\theta}+c_\lambda h_{\lambda,\theta}$,
following the right-hand side of \eqref{eq:resolvent-representation},
with $u(v)=c_\lambda$. By Lemmas~\ref{gen:grad}--\ref{gen:dir}, $u$ is jointly Borel, uniformly Lipschitz along the rays and uniformly continuous at $v$ in the tree metric, and vanishes uniformly at radial infinity. Thus $u\in\mathscr B$. Moreover
\[
\mathcal{L}_\theta u=\lambda u-g
\]
for $r>0$, which extends to
\[
F:=\lambda u-g\in\mathscr B
\]
with $F(v)=\lambda c_\lambda-g(v)$. The definition of $c_\lambda$ gives
\[
\int_0^{2\pi}u_\theta'(0+)\nu(d\theta)=\eta(\lambda c_\lambda-g(v))=\eta F(v).
\]
Together with the preceding properties, the uniform derivative bounds imply that $u\in\mathcal D_W(\mathcal L)$. Proposition~\ref{prop:sticky-martingale-problem} therefore implies that $u(Y_t)-u(Y_0)-\int_0^tF(Y_s)ds$ is a local martingale. The stopped local martingale is a true martingale on each bounded time interval because $u$ and $F$ are bounded. As $F=\lambda u-g$, integration by parts gives
\[
\mathbb E_x[e^{-\lambda t}u(Y_t)]-u(x)
=-\mathbb E_x\!\left[\int_0^t e^{-\lambda s}g(Y_s)\,ds\right].
\]
The first expectation tends to zero as $t\to\infty$, while the second converges by dominated convergence. Hence
\[
u(x)=\mathbb E_x\int_0^\infty e^{-\lambda t}g(Y_t)dt=\mathcal R_\lambda g(x).
\]
Consequently the explicit formula is the probabilistic resolvent; in particular
\[
\mathcal R_\lambda\mathscr B\subset\mathscr B,
\]
$\mathcal R_\lambda$ is positive and $\|\mathcal R_\lambda\|\le1/\lambda$.

\textbf{Step 2 (resolvent identity and strong continuity).}
We first derive the resolvent identity without invoking the strong
Markov property after time change. For any
$f\in\mathcal D_W(\mathcal L)\cap\mathscr B$ with
$\mathcal Lf\in\mathscr B$, the martingale problem and the
boundedness of $f$ and $\mathcal Lf$ give, after exponential
discounting and letting the terminal time tend to infinity,
\[
f=\mathcal R_\lambda(\lambda f-\mathcal Lf),\qquad\lambda>0.
\]
Fix $g\in\mathscr B$ and $\lambda,\mu>0$. Step~1 shows that
$u=\mathcal R_\lambda g$ and $v=\mathcal R_\mu g$ belong to
$\mathcal D_W(\mathcal L)$ and satisfy
\[
(\lambda-\mathcal L)u=g,
\qquad
(\mu-\mathcal L)v=g.
\]
Thus $f=u-v\in\mathcal D_W(\mathcal L)$, and
$(\lambda-\mathcal L)f=(\mu-\lambda)v$.
Applying the discounted identity to $f$ yields the resolvent identity
\begin{equation}\label{eq:resolvent-identity}
\mathcal R_\lambda-\mathcal R_\mu=(\mu-\lambda)\mathcal R_\lambda \mathcal R_\mu
\qquad\text{on }\mathscr B.
\end{equation}
Let $\omega_g$ be the tree-uniform modulus of continuity of $g$. For any $\delta>0$, Lemma~\ref{gen:usc} and Markov's inequality imply
\[
\|\lambda \mathcal R_\lambda g-g\|_\infty\le \omega_g(\delta)+\frac{2\|g\|_\infty C_0}{\delta}\left(\frac{\Gamma(3/2)}{\sqrt\lambda}+\frac1\lambda\right).
\]
Hence
\[
\lambda \mathcal R_\lambda g\to g\quad\text{in }\mathscr B\text{ as }\lambda\to\infty.
\] We justify the injectivity and density needed for the resolvent form of the Hille--Yosida theorem. The resolvent identity \eqref{eq:resolvent-identity} implies that the operators $\mathcal R_\lambda$ commute and have the same range for all $\lambda>0$. Since $\mu \mathcal R_\mu g\to g$ for every $g\in\mathscr B$, this common range is dense in $\mathscr B$. Moreover, if $\mathcal R_\lambda g=0$ for some $\lambda>0$, the resolvent identity and commutativity give $\mathcal R_\mu g=0$ for every $\mu>0$, whence
\[
g=\lim_{\mu\to\infty}\mu \mathcal R_\mu g=0.
\]
Thus each $\mathcal R_\lambda$ is injective. Together with $\|\mathcal R_\lambda\|\le1/\lambda$, these facts allow us to apply the resolvent form of the Hille--Yosida theorem, which yields a unique strongly continuous contraction semigroup $(S_t)$ on $\mathscr B$ with resolvent $\mathcal R_\lambda$. For each $x$ and $g\in\mathscr B$, the scalar functions $t\mapsto (S_tg)(x)$ and $t\mapsto\mathbb E_xg(Y_t)$ are bounded, continuous in $t$ (the latter by path continuity), and have the same Laplace transform. Their equality follows from uniqueness of the Laplace transform. Thus $S_t=P_t$ on $\mathscr B$, proving invariance and strong continuity.

\textbf{Step 3 (both domain inclusions).}
The standard resolvent identity for a $C_0$-semigroup gives
\[
\begin{aligned}
D(\mathcal L_{\mathscr B})&=\operatorname{Ran}(\mathcal R_\lambda),\\
\mathcal L_{\mathscr B}\mathcal R_\lambda g&=\lambda \mathcal R_\lambda g-g.
\end{aligned}
\] By Step 1 and the lemmas, every $f\in D(\mathcal L_{\mathscr B})$ satisfies all the stated analytic and boundary conditions with $F=\mathcal L_{\mathscr B}f$.

Conversely, suppose $f\in\mathscr B$, each $f_\theta\in C^2((0,\infty))$, the raywise $\mathcal{L}_\theta f$ extends to $F\in\mathscr B$, and the integral Wentzell condition holds. Lemma~\ref{gen:grad} applied to
\[
(\lambda-\mathcal{L}_\theta)f_\theta=\lambda f_\theta-F_\theta
\]
yields uniform first/second derivative bounds and the existence of all $f'_\theta(0+)$. Therefore $f\in\mathcal D_W(\mathcal L)$, as defined in Definition~\ref{def:wentzell-domain}. By Proposition~\ref{prop:sticky-martingale-problem}, exponential discounting and passage to $t\to\infty$ give
\[
f=\mathcal R_\lambda(\lambda f-F).
\] Since $\lambda f-F\in\mathscr B$, we obtain $f\in D(\mathcal L_{\mathscr B})$ and $\mathcal L_{\mathscr B}f=F$. This completes the characterization.
\end{proof}

\begin{remark}[Scope of the generator characterization]
\label{rem:generator-scope}
Theorem~\ref{gen:gen} applies only under the additional
boundedness, uniform ellipticity, and radial Lipschitz
assumptions of this subsection. The martingale problem of
Proposition~\ref{prop:sticky-martingale-problem} and the
local Dynkin formula remain valid under the more general
hypotheses of the preceding sections.
\end{remark}

\section{Local properties}
\label{sec:local-properties}

We study small-radius exits from the vertex for sticky Walsh
diffusions with a general, possibly non-atomic, spinning measure and
ray-dependent coefficients. Using the local Wentzell boundary value
problem and Dynkin's formula, we first obtain an asymptotic estimate
for the mean exit time. We then derive an exact joint Laplace
transform of the exit time, vertex occupation time, and outgoing
direction. This yields an exponential occupation-time law, an explicit
angular exit distribution with total-variation convergence to the
spinning measure, and a joint small-radius limit. The exit-time and
occupation-time estimates extend related results of Berry and
Colantoni \cite[Proposition~3.2 and Lemmas~3.3--3.5]{BerryColantoni2026}
from finite star graphs to the present Walsh setting.

\subsection{Mean exit time from the vertex}
\label{subsec:mean-exit-time}

For $\delta>0$, recall the first hitting time of the radial
level $\delta$, defined by
\[
\tau_\delta
:=
\inf\{0\le t<\tau^Y:R_t^Y=\delta\}.
\]

Our first result describes the small-radius asymptotics of the mean
exit time when the process starts from the vertex.

\begin{prop}
\label{prop:mean-exit-time}
Assume that, for every $\theta\in[0,2\pi)$, the functions
$r\mapsto b_\theta(r)$ and $r\mapsto\sigma_\theta(r)$ are continuous
in a neighborhood of the vertex, and that
$b_\theta$ and $\sigma_\theta^{-1}$
are uniformly bounded there. Then, as $\delta\downarrow0$,
\begin{equation}\label{eq:mean-exit-asymptotic}
\mathbb E_v[\tau_\delta]
=
\eta\delta+O(\delta^2).
\end{equation}
In particular,
\[
\lim_{\delta\downarrow0}
\frac{\mathbb E_v[\tau_\delta]}{\delta}
=
\eta.
\]
Moreover, for every sufficiently small $\delta>0$,
\[
\mathbb P_v(\tau_\delta<\infty)=1,
\]
so the vertex $v$ is not a trap.
\end{prop}

\begin{proof}
Apply Proposition~\ref{prop:local-Wentzell-BVP} with
\(f_\theta(r)\equiv -1.\)
Let $u=(u_\theta)_{\theta\in [0, 2\pi)}$ be the corresponding solution.
Then
\[
\mathcal{L}_\theta u_\theta(r)=-1,
\qquad 0<r<\delta,
\]
and
\(u_\theta(\delta)=0.\)
Moreover, the values $u_\theta(0)$ are independent of $\theta$; we
denote their common value by $u(v)$. Since the forcing term is
constant and equal to $-1$, its common value at the vertex is
\(\mathcal{L}u(v)=-1.\)
The Wentzell boundary condition therefore takes the form
\[
\int_0^{2\pi}u_\theta'(0+)\,\nu(d\theta)
=
-\eta.
\]

The solution $u$ of the local Wentzell boundary-value problem
satisfies the hypotheses of the local Dynkin formula,
Proposition~\ref{prop:local-Dynkin}. Indeed, the estimates of
Proposition~\ref{prop:local-Wentzell-BVP} imply the required
uniform bounds on $u$ and its radial derivatives on
$\mathcal W_\delta$, while
\(\mathcal{L}_\theta u_\theta=-1,
\qquad 0<r<\delta,\)
and the Wentzell boundary condition is satisfied at the vertex.

We first show that
\(\mathbb E_v[\tau_\delta]<\infty\)
without assuming this fact a priori.
By Proposition~\ref{prop:local-Dynkin},
the stopped Dynkin identity \eqref{eq:stopped-Dynkin}
applies to $u$. Since $\mathcal Lu=-1$ on
$\mathring{\mathcal W}_\delta$ (including at the vertex),
we obtain, for every $t\ge0$,
\[
\mathbb E_v
\bigl[u(Y_{t\wedge\tau_\delta})\bigr]
=
u(v)
-
\mathbb E_v[t\wedge\tau_\delta].
\]
Therefore,
\[
\mathbb E_v[t\wedge\tau_\delta]
=
u(v)
-
\mathbb E_v
\bigl[u(Y_{t\wedge\tau_\delta})\bigr].
\]
Since \(u\) is bounded on \(\mathcal W_\delta\),
\[
\mathbb E_v[t\wedge\tau_\delta]
\le
|u(v)|+\|u\|_{\infty,\delta},
\]
uniformly in \(t\). Letting \(t\to\infty\) and using monotone
convergence gives
\(\mathbb E_v[\tau_\delta]<\infty.\) We may therefore apply
the terminal-time identity \eqref{eq:terminal-Dynkin} to $u$.
Since
\(u_\theta(\delta)=0,
\; \theta\in[0,2\pi),\)
we have
\(u(Y_{\tau_\delta})=0.\)
Moreover,
\(\mathcal{L}u=-1\,\text{on }\mathring{\mathcal W}_\delta.\)
Hence
\[
0
=
u(v)-\mathbb E_v[\tau_\delta],
\]
and therefore
\(\mathbb E_v[\tau_\delta]=u(v).\)
It remains to determine the small-$\delta$ behavior of $u(v)$.
By Proposition~\ref{prop:local-Wentzell-BVP},
\[
\sup_{\theta\in [0, 2\pi)}
\sup_{0<r<\delta}
|u_\theta''(r)|
=
O(1)
\qquad\text{as }\delta\downarrow0.
\]
Taylor's formula at the vertex therefore gives, uniformly in
$\theta$,
\[
u_\theta(\delta)
=
u(v)
+
\delta u_\theta'(0+)
+
O(\delta^2).
\]
Since $u_\theta(\delta)=0$,
\[
0
=
u(v)
+
\delta u_\theta'(0+)
+
O(\delta^2).
\]
Integrating with respect to $\nu(d\theta)$ yields
\[
0
=
u(v)
+
\delta
\int_0^{2\pi}
u_\theta'(0+)\,\nu(d\theta)
+
O(\delta^2).
\]
Using the Wentzell boundary condition
established above, we obtain
\(u(v)=\eta\delta+O(\delta^2).\) Finally, since $\mathbb E_v[\tau_\delta]=u(v)$,
we obtain \eqref{eq:mean-exit-asymptotic}.
The stated limit follows immediately.
\end{proof}

\subsection{Joint exit-time and vertex occupation-time distributions}
\label{subsec:joint-exit-occupation}

We complement the preceding small-radius estimates with an exact
joint Laplace transform of the exit time, the time spent at the
vertex, and the outgoing direction. Related exit, local-time, and occupation-time questions have
been studied for one-dimensional diffusions and diffusions on graphs;
see, for instance, \cite{PitmanYor2003,BenichouDesbois2009}.
Salminen and Stenlund \cite{SalminenStenlund2025} also consider
angular-sector occupation times for Walsh Brownian motion with
a continuous spinning distribution. Related joint distribution
results for one-dimensional skew-sticky and sticky Brownian motion
appear in \cite{Touhami2021,CasterasMonsaingeon2026}.
The transform below instead concerns the first radial exit,
vertex occupation time, and exit direction for sticky Walsh
diffusions with arbitrary spinning measures and ray-dependent
coefficients.

For a process started from $v$, put
\[
\Gamma_\delta
:=\int_0^{\tau_\delta}\mathbf1_{\{Y_s=v\}}\,ds.
\]
For $\lambda\ge0$ and $\theta\in[0,2\pi)$, let
$\phi_{\lambda,\theta}$ and $\psi_{\lambda,\theta}$ be the
solutions on a fixed sufficiently small neighborhood of $0$ to
\[
\mathcal{L}_\theta w=\lambda w,
\]
with the initial conditions
\[
\begin{aligned}
\phi_{\lambda,\theta}(0)&=1,&
\phi'_{\lambda,\theta}(0)&=0,\\
\psi_{\lambda,\theta}(0)&=0,&
\psi'_{\lambda,\theta}(0)&=1.
\end{aligned}
\]

\begin{prop}
\label{prop:joint-exit-transform}
Let $Y$ be a sticky Walsh diffusion associated with
$(b,\sigma,\nu,\eta)$, started from $v$. Assume the local
continuity and uniform boundedness conditions of
Proposition~\ref{prop:mean-exit-time}. For sufficiently small
$\delta>0$, define
\[
A_\lambda(g;\delta)
:=\int_0^{2\pi}
\frac{g(\theta)}{\psi_{\lambda,\theta}(\delta)}
\,\nu(d\theta),
\qquad
B_\lambda(\delta)
:=\int_0^{2\pi}
\frac{\phi_{\lambda,\theta}(\delta)}
{\psi_{\lambda,\theta}(\delta)}\,\nu(d\theta).
\]
Then, for every $\lambda,q\ge0$ and every bounded Borel function
$g:[0,2\pi)\to\mathbb R$,
\begin{equation}\label{eq:joint-exit-transform}
\mathbb E_v\!\left[
 e^{-\lambda\tau_\delta-q\Gamma_\delta}
 g(\Theta_{\tau_\delta})\right]
=
\frac{A_\lambda(g;\delta)}
{B_\lambda(\delta)+\eta(\lambda+q)}.
\end{equation}
\end{prop}

\begin{proof}
Fix $\lambda,q\ge0$ and $0<\delta$ sufficiently small. By
joint Borel measurability and the uniform local bounds on the
coefficients, the solutions of the initial-value problems are
jointly Borel measurable in $(r,\theta)$. Their integral-equation
representations yield, uniformly in $\theta$,
\[
\phi_{\lambda,\theta}(r)=1+O(r^2),
\qquad
\psi_{\lambda,\theta}(r)=r+O(r^2),
\qquad r\downarrow0.
\]
To justify positivity uniformly over all $\lambda\ge0$,
let $Q_\theta(r):=\int_0^r 2b_\theta(z)/\sigma_\theta^2(z)\,dz$.
Every solution of $\mathcal{L}_\theta w=\lambda w$ satisfies
\[
\bigl(e^{Q_\theta(r)}w'(r)\bigr)'
=\frac{2\lambda}{\sigma_\theta^2(r)}
e^{Q_\theta(r)}w(r).
\]
A first-zero argument therefore gives
$\phi_{\lambda,\theta}(r)\ge1$ and
$\psi_{\lambda,\theta}(r)\ge s_\theta(r)>0$
for $r>0$, where $s_\theta(r):=\int_0^r e^{-Q_\theta(y)}\,dy$.
In particular, $\psi_{\lambda,\theta}(\delta)$ is
bounded away from zero uniformly in $\theta$, for each
$\delta>0$. Thus the integrals in the statement are finite
and $B_\lambda(\delta)>0$.

Write
\[
c:=\frac{A_\lambda(g;\delta)}
{B_\lambda(\delta)+\eta(\lambda+q)},
\qquad
a_\theta:=
\frac{g(\theta)-c\phi_{\lambda,\theta}(\delta)}
{\psi_{\lambda,\theta}(\delta)},
\]
and define, for $0\le r\le\delta$,
\[
u_\theta(r)
:=c\phi_{\lambda,\theta}(r)
+a_\theta\psi_{\lambda,\theta}(r).
\]
Then $u(v)=c$ is the common vertex value and
\[
\mathcal{L}_\theta u_\theta=\lambda u_\theta,
\qquad
u_\theta(\delta)=g(\theta),
\qquad
u_\theta'(0+)=a_\theta.
\]
The definition of $c$ gives
\begin{equation}\label{eq:discounted-vertex-condition}
\int_0^{2\pi}u_\theta'(0+)\,\nu(d\theta)
=A_\lambda(g;\delta)-cB_\lambda(\delta)
=\eta(\lambda+q)c.
\end{equation}
For fixed $\lambda$ and $\delta$, the integral-equation
representations and the uniform coefficient bounds control the
fundamental solutions and their first two radial derivatives
uniformly in $\theta$ on $[0,\delta]$. Since
$\inf_\theta\psi_{\lambda,\theta}(\delta)>0$ and $g$ is
bounded, the coefficients $a_\theta$, and hence $u$ and its first
two radial derivatives, are uniformly bounded there.

We verify explicitly the applicability of the generalized
Freidlin--Sheu formula. Choose a smooth radial cutoff
$\chi\in C^\infty([0,\infty))$ equal to one on $[0,\delta]$
and zero on $[2\delta,\infty)$, with $2\delta$ inside the
neighborhood on which the local coefficient bounds hold. Define
\[
\widetilde u_\theta(r)
=\begin{cases}
 u_\theta(r),&0\le r\le\delta,\\[1mm]
 \chi(r)\bigl[u_\theta(\delta)
 +(r-\delta)u_\theta'(\delta-)\bigr],&r>\delta.
\end{cases}
\]
The values and the first radial derivatives agree at $r=\delta$.
Thus $\widetilde u_\theta$ has locally absolutely continuous first
radial derivative, including across $r=\delta$, with uniformly
bounded second derivative almost everywhere on compact radial
intervals. The extension is bounded and has compact radial support;
its first derivatives at the vertex are uniformly bounded.
Furthermore, the joint Borel measurability of $u$, together with
that of $\theta\mapsto u_\theta'(\delta-)$, implies joint Borel
measurability of $\widetilde u$. The uniform first-derivative bound
also gives continuity at the common vertex in the tree topology.
Consequently $\widetilde u$ satisfies the radial test-function
conditions in \cite[Definition~2.6]{KaratzasYan2019}. The
generalized Freidlin--Sheu formula
\cite[Theorem~2.12]{KaratzasYan2019} can therefore be applied
directly to the sticky process $Y$, stopped at $\tau_\delta$:
the radial semimartingale equation and the local-time partition
property required for this application are part of
Definition~\ref{def:sticky-Walsh-diffusion}.
Although \eqref{eq:discounted-vertex-condition} differs, when
$q>0$, from the usual Wentzell condition for an eigenfunction,
the Freidlin--Sheu formula itself does not require that boundary
condition. The extension agrees with $u$ up to radial exit.

For $t\ge0$, put $s_t:=t\wedge\tau_\delta$ and
$H_t:=\int_0^{s_t}\mathbf1_{\{Y_z=v\}}\,dz$.
The generalized Freidlin--Sheu formula and
\eqref{eq:discounted-vertex-condition} yield
\[
u(Y_{s_t})-c
=N_{s_t}
+\lambda\int_0^{s_t}\mathbf1_{\{Y_z\neq v\}}u(Y_z)\,dz
+\eta(\lambda+q)c\,L_{s_t}^0(R^Y),
\]
where $N$ is a continuous local martingale. By the sticky
occupation identity, $\eta L_{s_t}^0(R^Y)=H_t$.
Since $u(Y_z)=c$ on $\{Y_z=v\}$, this becomes
\[
u(Y_{s_t})-c
=N_{s_t}
+\lambda\int_0^{s_t}u(Y_z)\,dz
+qcH_t.
\]
Integration by parts now shows that
\[
Z_t:=e^{-\lambda s_t-qH_t}u(Y_{s_t})
\]
is a continuous local martingale. Because $u$ is bounded on
$\mathcal W_\delta$ and the discount factor is at most one,
$Z$ is bounded and hence a true martingale.

By Proposition~\ref{prop:mean-exit-time},
$\mathbb E_v[\tau_\delta]<\infty$, so
$\tau_\delta<\infty$ almost surely. Letting $t\to\infty$
and using bounded convergence, $H_t\to\Gamma_\delta$ and
$u(Y_{\tau_\delta})=g(\Theta_{\tau_\delta})$ give
\[
c=\mathbb E_v\!\left[
 e^{-\lambda\tau_\delta-q\Gamma_\delta}
 g(\Theta_{\tau_\delta})\right].
\]
The definition of $c$ proves \eqref{eq:joint-exit-transform}.
\end{proof}

\begin{corollary}
\label{cor:joint-small-exit}
Under the assumptions of
Proposition~\ref{prop:joint-exit-transform}, let
\[
s_\theta(r):=\int_0^r
\exp\!\left(-2\int_0^y
\frac{b_\theta(z)}{\sigma_\theta^2(z)}\,dz\right)dy,
\qquad
K_\delta:=\int_0^{2\pi}
\frac{1}{s_\theta(\delta)}\,\nu(d\theta).
\]
For every sufficiently small $\delta>0$:
\begin{enumerate}
\item The vertex occupation time has the exponential law
\[
\Gamma_\delta\sim
\operatorname{Exp}\!\left(\frac{K_\delta}{\eta}\right),
\qquad
\mathbb E_v[\Gamma_\delta]=\frac{\eta}{K_\delta},
\]
where $\operatorname{Exp}$ is parametrized by its rate.
In particular,
\[
\mathbb E_v[\Gamma_\delta]
=\eta\delta+O(\delta^2),\qquad \delta\downarrow0.
\]
\item The random variables $\Gamma_\delta$ and
$\Theta_{\tau_\delta}$ are independent, and
\[
\mathbb P_v(\Theta_{\tau_\delta}\in d\theta)
=\frac{s_\theta(\delta)^{-1}}{K_\delta}\,\nu(d\theta).
\]
Moreover, as $\delta\downarrow0$,
\[
\left\|\mathbb P_v(\Theta_{\tau_\delta}\in\cdot)-\nu\right\|_{\mathrm{TV}}
=O(\delta).
\]
\item If $E\sim\operatorname{Exp}(1)$ and
$\Theta$ has distribution $\nu$, with $E$ and $\Theta$
independent, then
\[
\left(\frac{\tau_\delta}{\delta},
\Theta_{\tau_\delta}\right)
\xrightarrow[\delta\downarrow0]{\mathcal{L}}
(\eta E,\Theta)
\qquad\text{on }[0,\infty)\times S^1.
\]
\end{enumerate}
\end{corollary}

\begin{proof}
For $\lambda=0$, the initial-value problems give
$\phi_{0,\theta}\equiv1$ and $\psi_{0,\theta}=s_\theta$.
The transform identity \eqref{eq:joint-exit-transform} yields
\[
\mathbb E_v\!\left[
 e^{-q\Gamma_\delta}g(\Theta_{\tau_\delta})\right]
=
\frac{\displaystyle\int_0^{2\pi}
 g(\theta)s_\theta(\delta)^{-1}\,\nu(d\theta)}
{K_\delta+\eta q}.
\]
Taking $g\equiv1$ identifies the exponential law. Taking
$q=0$ gives the angular distribution, and the factorization
of the right-hand side gives independence.

Uniform boundedness of $b_\theta/\sigma_\theta^2$ near zero
gives, uniformly in $\theta$,
\[
s_\theta(\delta)=\delta+O(\delta^2),
\qquad
K_\delta=\frac{1}{\delta}+O(1).
\]
Consequently,
\[
\sup_{\theta\in[0,2\pi)}
\left|\frac{s_\theta(\delta)^{-1}}{K_\delta}-1\right|
=O(\delta),
\]
which proves the total-variation estimate in part~(2).
It also follows that $\Gamma_\delta/\delta$ converges in law to
$\eta E$. By independence, the pair
$(\Gamma_\delta/\delta,\Theta_{\tau_\delta})$
converges jointly to $(\eta E,\Theta)$.
Finally, Proposition~\ref{prop:mean-exit-time} and the exact
occupation-time mean established in part~(1) give
\[
0\leq \mathbb E_v[\tau_\delta-\Gamma_\delta]=\mathbb E_v[\tau_\delta]-\frac{\eta}{K_\delta}=O(\delta^2).
\]
As the difference is nonnegative, it follows that
$(\tau_\delta-\Gamma_\delta)/\delta\to0$ in $L^1$.
The claimed joint limit follows by Slutsky's theorem.
\end{proof}

\begin{remark}
The exact angular exit distribution has the same inverse-scale
weighting as in the non-sticky Walsh setting; compare
\cite[Theorem~4.5(i)]{KaratzasYan2019}.
In particular, this distribution does not depend on the
stickiness parameter $\eta$, whereas the exponential law
of $\Gamma_\delta$ has rate $K_\delta/\eta$.
The independence in Corollary~\ref{cor:joint-small-exit}
concerns $\Gamma_\delta$ and the exit direction;
it does not imply independence of the total exit time
$\tau_\delta$ and the exit direction when the radial
coefficients depend on $\theta$.
\end{remark}

\end{document}